\documentclass[aos,MSNbibl,citesort,dvips]{arximspdf}

\doi{10.1214/11-AOS924}
\volume{39}
\issue{6}
\pubyear{2011}
\firstpage{2883}
\lastpage{2911}

\makeatletter

\newtheorem{theorem}{Theorem}

\newproclaim{condition}{Condition}
\newproclaim{definition}{Definition}

\newtheorem{lemma}{Lemma}

\newtheorem{proposition}{Proposition}

\newproclaim{remark}{Remark}

\makeatother

\begin{document}
\begin{frontmatter}

\title{Rates of contraction for posterior distributions in~$\bolds{L^r}$-metrics, $\bolds{1 \le r \le\infty}$}
\runtitle{$L^r$ and uniform consistency of Bayes estimates}

\begin{aug}
\author[A]{\fnms{Evarist} \snm{Gin\'{e}}\corref{}\ead[label=e1]{gine@math.uconn.edu}}
\and
\author[B]{\fnms{Richard} \snm{Nickl}\ead[label=e2]{r.nickl@statslab.cam.ac.uk}}
\runauthor{E. Gin\'{e} and R. Nickl}
\affiliation{University of Connecticut and University of Cambridge}
\address[A]{Department of Mathematics \\
University of Connecticut \\
Storrs, Connecticut 06269-3009\\
USA\\
\printead{e1}}
\address[B]{Statistical Laboratory \\
Department of Pure Mathematics \\
\quad and Mathematical Statistics \\
University of Cambridge \\
Wilberforce Road\\
CB3 0WB, Cambridge\\
United Kingdom \\
\printead{e2}} 
\end{aug}

\received{\smonth{3} \syear{2011}}
\revised{\smonth{9} \syear{2011}}

%
\begin{abstract}
The frequentist behavior of nonparametric Bayes estimates, more
specifically, rates of contraction of the posterior distributions to
shrinking $L^r$-norm neighborhoods, $1 \le r \le\infty$, of the
unknown parameter, are studied. A theorem for nonparametric density
estimation is proved under general approximation-theoretic assumptions
on the prior. The result is applied to a variety of common examples,
including Gaussian process, wavelet series, normal mixture and
histogram priors. The rates of contraction are minimax-optimal for $1
\le r \le2$, but deteriorate as $r$ increases beyond $2$. In the case
of Gaussian nonparametric regression a Gaussian prior is devised for
which the posterior contracts at the optimal rate in all $L^r$-norms,
$1 \le r \le\infty$.
\end{abstract}

%

\end{frontmatter}

\section{Introduction}

In finite-dimensional statistical models the Bernstein--von Mises
theorem provides a frequentist justification of the use of Bayesian
methods. In the case of infinite-dimensional models, consistency
properties in weak metrics hold under relatively mild conditions; see
Schwartz \cite{S66}. Consistency in stronger metrics was considered by
Barron, Schervish and Wasserman \cite{BSW99} and by Ghosal, Ghosh and
Ramamoorthi \cite{GGR99}, and, shortly after, Ghosal, Ghosh and van der
Vaart \cite{GGV00} and Shen and Wasserman \cite{SW01} developed
techniques that allow us to prove frequentist rates of contraction of
the posterior to the true infinite-dimensional parameter in the
Hellinger metric, if the prior is suitably chosen according to the
structure of the nonparametric problem at hand. This led to further
progress recently; we refer to \cite{GV07,GV07b,VV08,VV09} and the
references therein.

This literature has been successful in generalizing the scope of these
techniques to a variety of different statistical models, and has
naturally focussed on consistency and rates of contraction results in
the \textit{Hellinger} distance. For instance, if $p_0$ is the unknown
density to be estimated, and if $\Pi(\cdot|X_1,\ldots,X_n)$ is the
posterior based on a prior $\Pi$ and a sample $X_1,\ldots, X_n$ with
joint law~$P^n_0$, results of the kind
%
\begin{equation} \label{cont}
\Pi\bigl(p\dvtx h(p,p_0) \ge\varepsilon_n |X_1,\ldots, X_n\bigr) \to
0 \qquad\mbox{in } P^n_0 \mbox{ probability}
\end{equation}
were established, where $h^2(f,g) = \int(\sqrt f - \sqrt g)^2$ is the
Hellinger metric and where $\varepsilon_n\to0$. Such posterior
contraction results are known to imply the same frequentist consistency
rate $\varepsilon_n$, also in the metric $h$, for the associated formal
Bayes estimators.

In this article we investigate the question of how to generalize
results of this kind to more general loss-functions than the Hellinger
metric, with a~particular focus on $L^r$-norms, $1 \le r \le\infty$.
Such results are of interest for a~variety of reasons, for example, the
construction of simultaneous confidence bands, or for plug-in
procedures that require control of nonparametric remainder terms (e.g.,
in the proof of the Bernstein--von Mises theorem in semiparametric
models in Castillo \cite{C11}). They are also of interest with a view
on a more unified understanding of nonparametric Bayes procedures that
complements the existing $L^r$-type results for standard frequentist methods.

The main challenge in extending the theory to the $L^r$-case, except
for specific conjugate situations discussed below, rests in
generalizing the Le Cam--Birg\'{e} testing theory for the Hellinger
metric to more general situations. A main ingredient of the proof of a
result of the kind (\ref{cont}) is that, in testing problems of the form
%
\begin{equation} \label{bltest}
H_0\dvtx p=p_0  \quad\mbox{against}\quad H_A\dvtx p \in\{p\dvtx h(p,p_0) \ge\varepsilon
_n\},
\end{equation}
universal tests with concentration bounds on type-II errors of the type~$e^{-Cn\varepsilon_n^2}$ exist, under assumptions on the size, or
entropy, of the ``alternative'' space defining~$H_A$. This fact is
rooted in the subtle connection between nonparametric testing problems
and the Hellinger metric as highlighted in the work of Le Cam \cite
{LC86} and Birg\'{e}~\cite{B83}. A main contribution of this article is
the development of a new approach to testing problems of the kind (\ref
{bltest}) based on concentration properties of linear centered
kernel-type density estimators, derived from empirical process
techniques. While this approach can only be used if one has sufficient
control of the approximation properties of the support of the prior, it
can be generalized to arbitrary $L^r$-metrics, including the supremum
norm $\|f\|_\infty= {\sup_x} |f(x)|$. The concentration properties of
these tests depend on the geometry of the $L^r$-norm and deteriorate as
$r \to\infty$, which is, in a sense, dual to the fact that the minimax
testing rate in the sense of Ingster \cite{I93} approaches the minimax
rate of estimation as $r \to\infty$.

While our main results can be viewed as ``abstract'' in that they
replace the entropy conditions in \cite{GGV00} for sieve sets $\mathcal
P_n$ by general approximation-theoretic conditions (see Theorems \ref
{general} and \ref{generalbd} below), our findings become most
transparent by considering specific examples, selected in an attempt to
reflect the spectrum of situations that can arise in Bayesian
nonparametrics: In Section \ref{ex} we study the ``ideal'' situation of
a simple uniform wavelet prior on a H\"{o}lder ball, the
``supersmooth''
situation of mixtures of normals, the case of random histograms based
on a Dirichlet process where no uniform bound on the $L^\infty$-norm of
the support of the prior is available, as well as Gaussian process
priors of the kind studied in~\cite{VV08}. The general conclusion is
that if $f_0$ is $\alpha$-smooth, then the rate of contraction obtained
in the $L^r$-norm for a posterior based on an adequately chosen prior
of smoothness $\alpha$ is, up to $\log n$ factors, and with $\bar r
=\max(2,r)$,
%
\begin{equation} \label{rate0}
\biggl(\frac{1}{n}\biggr)^{({\alpha- 1/2 + 1/\bar r})/({2\alpha+1})}.
\end{equation}
So as soon as $r \le2$ our proof retrieves the minimax optimal rate,
but for $r>2$ the rate deteriorates by a genuine power of $n$. As
$\alpha$ approaches infinity this effect becomes more lenient and
vanishes in the limit.

We currently have no proof of the fact that our general theorem gives
the right rate for Bayesian posteriors if $r>2$---similar problems are
known with nonparametric maximum likelihood estimators in $L^r$-metrics
(cf. the proof of Proposition~6 in \cite{N07}). While we do not settle
the issue of optimality of our rates for $r>2$ in this article, we also
prove in Theorem \ref{whitenoise1} below that in nonparametric Gaussian
regression the minimax rate of contraction can be obtained by certain
diagonal Gaussian wavelet priors, in all $L^r$-norms simultaneously. We
believe that this result is closely tied to the fact that the posterior
is then itself Gaussian, and conjecture that our rates cannot be
substantially improved in the nonconjugate situation.

\section{Main results} \label{ex}

Let $\mathcal P$ be a class of probability densities on $[0,1]$ or
$\mathbb R$, and let $X_1,\ldots,X_n$ be a random sample drawn from some
unknown probability density $p_0$ with joint law the first $n$
coordinate projections of the infinite product probability measure
$P_0^\mathbb N$. Suppose one is given a prior probability distribution
$\Pi$ defined on some $\sigma$-algebra $\mathcal B$ of $\mathcal P$.
The posterior is the random probability measure
\[
\Pi(B|X_1, \ldots, X_n)= \frac{\int_B \prod_{i=1}^n p(X_i)\,d\Pi(p)}{\int
_\mathcal P \prod_{i=1}^n p(X_i)\,d\Pi(p)}, \qquad  B \in\mathcal B.
\]
We wish to analyze contraction properties of the posterior distribution
under certain regularity conditions on $\Pi$ and $p_0$, and these
regularity properties can be conveniently characterized by wavelet theory.

\subsection{Function spaces and wavelets}\label{w}

For $T=\mathbb R$ or $T=[0,1]$, $f\dvtx T\mapsto\mathbb R$, we shall write
$\|f\|_\infty={\sup_{x \in T}}|f(x)|$, the norm on the space $C(T)$ of
bounded continuous real-valued functions defined on $T$. We shall use
wavelet theory throughout; see\vadjust{\goodbreak} \cite{M92,HKPT}. Let $\phi, \psi$ be the
scaling function and wavelet of a~multiresolution analysis of the space
$L^2(T)$ of square integrable real-valued functions on $T$. We shall
say that the wavelet basis is $S$-\textit{regular} if $\phi, \psi$ are
$S$-times continuously differentiable on $T$. For instance we can take
Daubechies wavelets on $T=\mathbb R$ of sufficiently large order $N$
(see \cite{M92}) and define the translated scaling functions and wavelets
%
\begin{equation} \label{trans}
\phi_k = \phi(\cdot-k),\qquad \psi_{\ell k} = 2^{\ell/2}\psi\bigl(2^\ell(\cdot
)-k\bigr),\qquad \ell\in\mathbb N \cup\{0\}, k \in\mathbb Z,
\end{equation}
which form an orthonormal basis of $L^2(\mathbb R)$.

For $T=[0,1]$ we consider the orthonormal wavelet bases of
$L^2([0,1])$ constructed in Theorem 4.4 of Cohen, Daubechies and Vial
\cite{CDV93}. Each such basis is built from a Daubechies scaling
function $\phi$ and its corresponding wavelet~$\psi$, of order $N$,
starting at a fixed resolution level $J_0$ such that $2^{J_0}\ge2N$
(see Theorem~4.4 in \cite{CDV93}): the $\psi_{\ell k}, \phi_k$ that are
supported in the interior of $[0,1]$ are all kept, and suitable
boundary corrected wavelets are added, so that the $\{\phi_k, \psi
_{\ell k}\dvtx 0 \le k < 2^\ell, \ell \in\mathbb N, \ell > J_0\}$ still form an
orthonormal basis for~$L^2([0,1])$. While formula (\ref{trans}) now
only applies to the ``interior'' wavelets, one can still write $\phi
_{jk}= 2^{j/2} \phi_k(2^j \cdot)$ for every $k, j \ge J_0$; cf. page 73
in \cite{CDV93} and also after Condition \ref{kernel} below.
\begin{definition} \label{besov}
Let $T=[0,1]$ or $T=\mathbb R$, and let $1\leq p,q \leq\infty$, $0\le
s<S$, \mbox{$s \in\mathbb R$}, $S \in\mathbb N$. Let $\phi, \psi$ be
bounded, compactly supported $S$-regular scaling function and wavelet,
respectively, and denote by $\alpha_k(f)=\int_T \phi_k f$ and $\beta
_{\ell k}(f)=\int_T \psi_{\ell k}f$ the wavelet coefficients of $f \in
L^p(T)$. The Besov spa-\break ce~$B^s_{pq}(T)$ is defined as the set of
functions $\{f \in L^p(T)\dvtx \|f\|_{s,p,q} <\infty\}$ where
\[
\|f\|_{s,p,q} :=\bigl\|\alpha_{(\cdot)}(f)\bigr\|_p + \Biggl(\sum_{\ell=0}^\infty
\bigl(2^{\ell(s+1/2-1/p)}\bigl\|\beta_{\ell(\cdot)}(f)\bigr\|_p \bigr)^q \Biggr)^{1/q}
\]
with the obvious modification in case $q=\infty$.
\end{definition}
\begin{remark} \label{besovmania}
We note the following
standard embeddings/identifications we shall use (cf. \cite{M92,HKPT}):
for $\mathcal C^s(T)$ the H\"{o}lder (-Zygmund in case $s$ integer)
spaces on $T$, we have $B^s_{\infty\infty}(T)=\mathcal C^s(T)$.
Moreover $B^s_{22}(T) = H^s(T)$ where~$H^s(T)$ are the standard
$L^2$-Sobolev spaces. We also have the ``Sobolev-type'' imbeddings
$B^s_{rq}(T) \subset B_{tq}^{s-1/r+1/t}(T)$ for $t \ge r, 1\le q\le
\infty$. Finally,\vspace*{1pt} if $T=[0,1]$, then $C^\alpha(T) \subset B^\alpha_{r
\infty}(T)$ for every\vspace*{1pt} $r \le\infty$, where $C^\alpha(T)=\{f\dvtx T\mapsto
\mathbb R\dvtx \|f\|_{\alpha,\infty}<\infty\}$, with $\|f\|_{\alpha,\infty
}:={\sum_{k=0}^{\alpha}}\|f^{(k)}\|_\infty$, $\alpha\in\mathbb N$.
\end{remark}

\subsection{Uniform wavelet series}

Let us consider first the case where an a~priori upper bound on the H\"
{o}lder norm $\|p_0\|_{\alpha, \infty, \infty}$ is available, so that
the prior can be chosen to have bounded support in $\mathcal{C}^{\alpha
}([0,1])$. An example is obtained, for example, by uniformly
distributing wavelet coefficients on a~H\"{o}lder ball. Let $\{\phi
_k, \psi_{\ell k} \}$ be a $N$-regular CDV-wavelet basis for
$L^2([0,1])$, let $u_{\ell k}$ be i.i.d.~$U(-B,B)$ random variables,
and define, for $\alpha<N$, the random wavelet series
%
\begin{equation}
U_\alpha(x) = \sum_k u_{0k} \phi_k(x) +\sum_{\ell=J_0}^\infty\sum_k
2^{-\ell(\alpha+1/2)}u_{\ell k} \psi_{\ell k}(x),
\end{equation}
which has trajectories in $\mathcal C^{\alpha}([0,1]) \subset
L^r([0,1]), 1 \le r \le\infty$, almost surely (in view of Definition
\ref{besov} and Remark \ref{besovmania}). Since moreover $\|U_{\alpha}\|
_\infty\le C(B,\alpha, \psi)$, and since the exponential map has
bounded derivatives on bounded subsets of $\mathbb R$, the same applies
to the random density
\[
p^{U, \alpha}(x):=\frac{e^{U_\alpha(x)}}{\int_0^1 e^{U_{\alpha}(y)}\,dy},
\]
whose induced law on $C([0,1])$ we denote by $\Pi^\alpha$. Our general
results below imply the following proposition,
which, since $p_0$ is bounded away from zero, implies the same
contraction rate in Hellinger distance $h$. Note moreover that the
result for $2<r<\infty$ could be obtained from interpolation properties
of $L^r$-spaces.
\begin{proposition} \label{unif}
Let $X_1, \ldots, X_n$ be i.i.d. on $[0,1]$ with density $p_0$
satisfying $\|{\log p_0}\|_{\alpha, \infty} \le B$. Let $1 \le r \le
\infty$, $\bar r = \max(2,r), r^* = \min(r,2)$, and suppose $\alpha\ge
1-1/r^*$. Then there exist finite positive constants $M, \eta=\eta
(\alpha,r)$ such that, as $n \to\infty$,
%
\begin{eqnarray}
&&\Pi^\alpha\bigl\{p\in\mathcal{P}\dvtx\|p-p_0\|_r\ge Mn^{-({\alpha
-1/2+1/\bar r})/({2\alpha+1})}(\log n)^\eta|X_1,\ldots,X_n\bigr\}\nonumber\\[-8pt]\\[-8pt]
&&\qquad\to
^{P_0^{\mathbb N}} 0.\nonumber
\end{eqnarray}
\end{proposition}

\subsection{Dirichlet mixtures}

Consider first, as in  \cite{GGR99,GV01,GV07b}, a normal mixture prior
$\Pi$, defined as follows: for $\varphi$ the standard normal density,
set:\vspace*{8pt}

(-) $p_{F,\sigma} = \int_\mathbb R \sigma^{-1} \varphi((\cdot-y)/\sigma)
\,dF(y)$,

(-) $F \sim D_\alpha$ the Dirichlet-process with base measure $\alpha=
\alpha(\mathbb R)\bar\alpha$, $\alpha(\mathbb R)<\infty$ and $\bar
\alpha$ a probability measure,

(-) $\sigma\sim G$, where $G$ is a probability distribution with
compact support in $(0, \infty)$.
\begin{proposition} \label{dir}
Let $X_1, \ldots, X_n$ be i.i.d. on $\mathbb R$ with density $p_{F_0,
\sigma_0}$ where $\sigma_0>0$ and where $F_0$ is supported in $[-k_0,
k_0], k_0>0$. Suppose that $G$ has a positive continuous density in a
neighborhood of $\sigma_0$, and that the base measure $\alpha$ has
compact support and a continuous density on an interval containing
$[-k_0, k_0]$. Then there exist finite positive constants $M, \eta$
such that
%
\begin{equation}
\Pi^\alpha\biggl\{p\in\mathcal{P}\dvtx\|p-p_0\|_\infty\ge M \frac{(\log
n)^\eta}{\sqrt n} \Big|X_1,\ldots,X_n\biggr\}\to^{P_0^{\mathbb N}} 0
\qquad\mbox{as }
n\to\infty.\hspace*{-28pt}
\end{equation}
\end{proposition}

Consider next a random histogram based on a Dirichlet process, similar
to the priors studied in \cite{S07}: for $j \in\mathbb N$ let $\mathrm{Dir}_j$
be a Dirichlet-distribution on the $2^j$-dimensional unit simplex, with
all parameters equal to one. Consider the dyadic random histogram with
resolution level $j$
\[
\sum_{k=1}^{2^j} \alpha_{jk} 2^{j}1\biggl\{\biggl(\frac{k-1}{2^j}, \frac
{k}{2^j} \biggr]\biggr\}(x),\qquad \{a_{jk}\} \sim \mathrm{Dir}_j,\qquad
x \in[0,1],
\]
and denote its law on the space of probability densities by $\Pi_j$.
Note that this prior is not concentrated uniformly (in $j$) on bounded
densities (despite the densities in the support being uniformly bounded
for fixed $j$).
\begin{proposition} \label{dir1}
Let $X_1, \ldots, X_n$ be i.i.d. on $[0,1]$ with density $p_0 \in
\mathcal C^{\alpha}([0,\allowbreak1]), 0<\alpha\le1$, satisfying\vspace*{1pt} $p_0>0$ on
$[0,1]$. Let $j_n$ be such that $2^{j_n} \sim(n/\allowbreak\log n)^{1/(2\alpha
+1)}$, let $1 \le r \le\infty$, $\bar r = \max(2,r)$ and let either
$\alpha>1/2$ or $r=1$. Then for some $M, \eta=\eta(\alpha,r)$, as $n\to
\infty$
%
\begin{eqnarray}
&&\Pi_{j_n}\bigl\{p\in\mathcal{P}\dvtx\|p-p_0\|_r \ge Mn^{-({\alpha
-1/2+1/\bar r})/({2\alpha+1})}(\log n)^\eta|X_1,\ldots,X_n\bigr\}\nonumber\\[-8pt]\\[-8pt]
&&\qquad\to
^{P_0^{\mathbb N}} 0.\nonumber
\end{eqnarray}
\end{proposition}

\subsection{Gaussian process priors} \label{gpp}

We now study a variety of Gaussian process priors that were considered
in the nonparametric Bayes literature recently; see \mbox{\cite{VV08,VV09}}
for references. To reduce technicalities we shall restrict ourselves to
integrated Brownian motions, but see also the remark below.
\begin{definition}\label{fbm}
Let $B(t)=B_{1/2}(t)$, $t\in[0,1]$, be a (sample-continuous version of)
standard Brownian motion. For $\alpha>1$, $\alpha\in\{n-1/2\dvtx n \in
\mathbb N\}$, setting $\{\alpha\}=\alpha-[\alpha]$, $[\alpha]$ being
the integer part of $\alpha$, $B_\alpha$ is defined as the $[\alpha
]$-fold integral
\begin{eqnarray*}
B_\alpha(t)&=& \int_0^t \int_0^{t_{[\alpha]-1}}\cdots\int_0^{t_2} \int
_0^{t_1}B(s)\,ds \,dt_1\cdots dt_{[\alpha]-1}\\
&=&\frac{1}{([\alpha]-1)!}\int_0^t(t-s)^{[\alpha]-1}B(s)\,ds,\qquad  t\in[0,1],
\end{eqnarray*}
where for $[\alpha]=1$ the multiple integral is understood to be only
$\int_0^tB(s)\,ds$.
\end{definition}

Following \cite{L91,VV08}, and as before Proposition \ref{unif}, we
would like to define our prior on densities as the probability law of
the random process
%
\begin{equation} \label{expo}
\frac{e^{B_\alpha}}{\int_0^1e^{B_\alpha(t)}\,dt},
\end{equation}
but we must make two corrections: first, since $B^{(k)}_\alpha(0)=0$
a.s., $k\le[\alpha]$, would impose unwanted\vadjust{\goodbreak} conditions on the value at
zero of the density, we should release $B_\alpha$ at zero, that is,
take $\bar B_\alpha:= \sum_{k=0}^{[\alpha]}Z_kt^k/k!+B_\alpha$,
where $Z_k$ are i.i.d.~$N(0,1)$ variables independent of $B_\alpha$;
see \cite{VV08}.
In order to deal with bounded densities, we introduce a second
modification to (\ref{expo}), and define our prior (on the Borel sets
of $C([0,1])$) as
%
\begin{equation}\label{prior}
\Pi=\mathcal{L}\biggl(\frac{e^{\bar{B}_\alpha}}{\int_0^1e^{\bar
{B}_\alpha(t)}\,dt}\Big| \|\bar{B}_\alpha\|_\infty\le c\biggr),
\end{equation}
where $c$ is a fixed arbitrary positive constant. This prior works as
follows: if $A\subset C([0,1])$ is a measurable set of continuous
densities on $[0,1]$, then
\[
\Pi(A)=\Pr\biggl\{ e^{\bar B_\alpha}\Big/\int e^{\bar B_\alpha}\in A, \| \bar
B_\alpha\|_\infty\le c\biggr\}\Big/\Pr\{ \| \bar B_\alpha\|_\infty\le c\},
\]
and clearly the denominator is strictly positive for all $c>0$; see
Proposition~\ref{smb} below.
\begin{proposition}\label{fbmp}
Let $1 \le r \le\infty$, $\bar r = \max(r,2), \alpha\in\{n-1/2, n
\in\mathbb N\}$ and assume \textup{(a)} $p_0\in\mathcal C^\alpha([0,1])$, and
\textup{(b)} $p_0$ is bounded and bounded away from zero, say, $2\|{\log p_0}\|
_\infty\le c<\infty$. Let $\Pi$ be the prior defined by (\ref{prior})
where $\alpha$ is as in \textup{(a)} and $c$ is as in \textup{(b)}. Then, if $X_i$ are
i.i.d. with common law $P_0$ of density $p_0$, there exists $M<\infty$ s.t.
\begin{eqnarray*}
&&\Pi\bigl\{p\in\mathcal{P}\dvtx\|p-p_0\|_r \ge Mn^{-({\alpha-1/2+1/\bar
r})/({2\alpha+1})}(\log
n)^{(1/2)1_{\{r=\infty\}}}|X_1,\ldots,X_n\bigr\}\\
&&\qquad\to0
\end{eqnarray*}
in $P_0^{\mathbb N}$-probability as $n\to\infty$.
\end{proposition}

As remarked before Proposition \ref{unif}, a contraction result in the
Hellinger distance follows as well, and the case $2<r<\infty$ could be
obtained from interpolation.

The result in Proposition \ref{fbmp} extrapolates to fractional
multiple integrals of Brownian motion (Riemann--Liouville processes) of
any real valued index $\alpha>1/2$, and it also extends to the related
fractional Brownian motion processes (see, e.g., \cite{VV08} for
definitions), but, for conciseness and clarity of exposition, we
refrain from carrying out these extensions.

\subsection{Sharp rates in the Gaussian conjugate situation}\label{conjug}

We currently have no proof that the rates obtained in the previous
subsections are optimal for these priors as soon as $r>2$. While we
conjecture that Bayesian posteriors may suffer from suboptimal
contraction rates in density estimation problems in $L^r$-loss, $r>2$,
we finally show here that in the much simpler conjugate situation of
nonparametric regression with Gaussian errors, sharp rates in all~$L^r$
norms can be obtained at least for certain diagonal wavelet priors. The
proof of this result follows from a direct analysis of the posterior
distribution, available in closed form due to conjugacy.\vadjust{\goodbreak}

Given a noise level $1/\sqrt n, n \in\mathbb N$, we observe
%
\begin{equation}\label{wnm}
dY^{(n)}(t)=f(t)\,dt+\frac{1}{\sqrt{n}}\,dB(t),\qquad t\in[0,1],
\end{equation}
for $f=f_0\in L^2([0,1])$, where $B$ is Brownian motion on $[0,1]$.
This model is well known to be asymptotically equivalent to
nonparametric regression with fixed, equally-spaced design and Gaussian errors.

Consider priors on $L^2([0,1])$ defined on a $S$-regular CDV-wavelet
basis as
%
\begin{equation}\label{finpri}
\Pi=\mathcal{L}\Biggl(\sum_{k=0}^N g_k\phi_k+\sum_{\ell=J_0}^\infty\sum
_{k=0}^{2^{\ell}-1}\sqrt{\mu_\ell}g_{\ell k}\psi_{\ell k}\Biggr)
\end{equation}
in $L^2([0,1])$, with the $g$'s i.i.d. $N(0,1)$ and with $\mu_\ell= \ell
^{-1} 2^{-\ell(2\alpha+1)}$ $\forall\ell\ge J_0$. Such a prior is
designed for $\alpha$-smooth $f_0$. As is easily seen, the series in
(\ref{finpri}) converges uniformly almost surely.
\begin{theorem}\label{whitenoise1}
Let $0<\alpha<S$, and let $\Pi$ be the Gaussian prior on $L^2([0,1])$
defined by (\ref{finpri}) based on a CDV wavelet basis of $L^2([0,1])$
of smoothness at least $S$. Let $f_0\in\mathcal C^\alpha([0,1])$, let
$\varepsilon_n = (n/\log n)^{-\alpha/(2 \alpha+1)}$ and suppose we
observe $dY_0^{(n)}(t)=f_0(t)\,dt+dB(t)/\sqrt{n}$. Then there exists
$C<\infty$ and $M_0<\infty$ depending only on the wavelet basis, $\alpha
$ and $\|f_0\|_{\alpha,\infty, \infty}$ such that, for every $M_0\le
M<\infty$, and for all $1 \le r \le\infty, n\in\mathbb N$,
%
\begin{equation}\label{contract2}
E_{Y_0^{(n)}}\Pi\bigl(f\dvtx\|f-f_0\|_r>M\varepsilon_n|Y_0^{(n)}
\bigr)\le n^{-C^2(M-M_0)^2}.
\end{equation}
\end{theorem}

This rate of convergence is sharp (in case $r<\infty$ up to the $\log
n$-term) in view of the usual minimax lower bounds and since the
contraction rate implies the same rate of convergence for the formal
Bayes estimator $E_\Pi(f|Y^{(n)}_0)$ to $f_0$ (using Anderson's lemma
and the fact that the posterior is a random Gaussian measure on
$L^2([0,1])$, as inspection of the proof shows). One may even apply the
usual thresholding techniques to the posterior mean to obtain a
Bayesian rate adaptive estimator of $f_0$ by proceeding as in \cite{GN09,LN10}.

\section{\texorpdfstring{General contraction theorems for density estimates in
$L^r$-loss, $1 \le r \le\infty$}
{General contraction theorems for density estimates in Lr-loss, 1 <= r <= infinity}}

We shall, in our main results, use properties of various approximation
schemes in function spaces, based on integrating a localized
kernel-type function $K_j(x,y)$ against functions $p$, $K_j(p)=\int
K_j(\cdot,y)p(y)\,dy$. Let, in slight abuse of notation, for $T \subseteq
\mathbb R$, $L^1(\mu_w)= L^1(T,\mathcal{B},\mu_w), w \ge0$ be the
space of $\mu_w$-integrable functions, $d\mu_w(t)=(1+|t|)^w\,dt$, normed
by $\|f\|_{\mu_w} = \int_T |f(t)| (1+|t|)^w\,dt$.
Recall the notion of $p$-variation of
a function (e.g., as before Lemma 1 in~\cite{GN09}).
\begin{condition} \label{kernel} Let $T=\mathbb R$ or $T=[0,1]$. The
sequence of operators $K_j(x,\allowbreak y)=2^jK(2^jx, 2^jy); x,y \in T, j
\ge0$,
is called an admissible approximating sequence if it satisfies one of
the following conditions:\vadjust{\goodbreak}

(a) (convolution kernel case): $K(x,y)=K(x-y)$, where $K \in L^\infty
(T)$ is of bounded $p$-variation for some finite $p \ge1$, right (or
left) continuous, and satisfies $\|K\|_{\mu_w}<\infty$ for some $w >2$.

(b) (multiresolution projection case): $K(x,y)=\sum_k \phi(x-k) \phi
(y-k)$, the sum extending over any subset of $\mathbb Z$, where $\phi
\in L^1 \cap L^\infty$ has bounded $p$-variation for some finite $p \ge
1$ and satisfies, in addition, ${\sup_{x \in\mathbb R} \sum_k }|\phi
_k(x)| < \infty$ as well as $|K(x,y)| \le\Phi(|x-y|)$ for every $x,y
\in T$ and some function $\Phi\in L^\infty(\mathbb R)$ for which
$\|\Phi\|_{\mu_w}<\infty$ for some $w >2$.

(c) (multiresolution case, $T=[0,1]$): $K(x,y)=\sum_k \phi_k(x) \phi
_k(y)$ is the projection kernel of a Cohen--Daubechies--Vial (CDV)
wavelet basis.
\end{condition}

Condition (a) is a standard assumption on kernels, condition (b) is
satisfied for most wavelet basis on $\mathbb R$, such as Daubechies,
Meyer or spline wavelets, by using standard wavelet theory (e.g., \cite
{HKPT}). For part (c) we note the following: as in the case of the whole
line, an orthonormal basis of $V_j=\{\phi_{jk} = 2^{j/2}\phi_k(2^j\cdot
)\}$ is obtained from $2^{j-J_0}$-fold dilates of the basic linear span
$V_{J_0}$, for every $j\ge J_0$ (pa\-ge~73 in \cite{CDV93}). In this
case, $V_j$ has dimension~$2^j$, and a basis consists of: (i)~$N$~left
edge functions $\phi_{jk}^0(x)=2^{j/2}\phi_k^0(2^jx)$, $k=0,\ldots,N-1$,
where $\phi_k^0$ is a modification of $\phi$, which is still bounded
and of bounded support; (ii) $N$ right edge functions $\phi
_{jk}^1(x)=2^{j/2}\phi_k^1(2^jx)$, $k=0,\ldots,N-1$, $\phi_k^1$ also
modifications\vspace*{-1pt} of $\phi$ bounded and of bounded support, and then the
$2^j-N$ ``interior'' usual translations of dilations of $\phi$, $\phi
_{jk}$, $k=N,\ldots, 2^j-N-1$.
The projection kernel $K_j(x,y)=K_j^0(x,y)+K_j^1(x,y)+\tilde K_j(x,y)$
corresponds to the projection onto the three orthogonal components of
$V_j$ (the linear spans, respectively, of the left edge functions $\phi
_{j,k}^0$, the right edge functions $\phi_k^1$, and the interior
functions $\phi_{jk}$). The first two spaces have dimension $N$ and the
third, $2^j-2N$. By Lemma 8.6 in~\cite{HKPT}, there exist bounded,
compactly supported nonnegative functions $\Phi$ such that $\tilde
K(x,y)\le\Phi(|x-y|)$, for all $x,y$. We call this function a
majorizing kernel of the interior part of $K$.

Let $X_i$ be i.i.d. with law $P_0$ and density $p_0$.
\begin{theorem}\label{general} Let $T=[0,1]$ or $T=\mathbb R$, let
$\mathcal{P}=\mathcal{P}(T)$ be a set of probability densities on $T$,
and let $\Pi_n$ be priors defined on some $\sigma$-algebra of $\mathcal
P$ for which the maps $p\mapsto p(x)$ are measurable for all $x\in T$.
Let $1\le r \le\infty$ and let $\varepsilon_n \to0$ as $n \to\infty$
be a sequence of positive numbers such that $\sqrt n \varepsilon_n \to
\infty$ as $n \to\infty$. Let
%
\begin{equation}\label{rate}
\delta_n = \varepsilon_n (n \varepsilon_n^2)^{{1}/{2} -
{1}/({2r})} \gamma_n
\end{equation}
for some sequence $\gamma_n$ satisfying $\gamma_n \ge1$ $\forall n$. Let
$J_n$ be any sequence satisfying $2^{J_n} \le c n \varepsilon_n^2$ for
some fixed $0<c<\infty$, and let $K_j$ be an admissible approximator
sequence. Let $\mathcal P_n$ be a sequence of subsets of
%
\begin{equation}\label{pr}
\{p \in\mathcal P\dvtx  \|K_{J_n}(p)-p\|_r\le C(K) \delta_n, \|p\|
_{\mu_w} \le D \},
\end{equation}
where $C(K)$ is a constant that depends only on the operator kernel
$K$, $D$ is a fixed constant, and where $w >(2-r)/r$ if $r<2$, $w=0$ if
$r \ge2$.

Assume there exists $C>0$ such that, for every $n$ large
enough:\vspace*{8pt}

\item{(1)} $\Pi_n(\mathcal{P}\setminus\mathcal{P}_n)\le
e^{-(C+4)n\varepsilon_n^2}$ and
\item{(2)} $\Pi_n\{p\in\mathcal{P}\dvtx -P_0\log\frac{p}{p_0}\le
\varepsilon_n^2, P_0(\log\frac{p}{p_0})^2\le\varepsilon
_n^2\}\ge e^{-Cn\varepsilon_n^2}$.\vspace*{8pt}

Let $p_0 \in L^r(T)$ be s.t. $\|K_{J_n}(p_0)-p_0\|_r =O(\delta_n)$ and
s.t. $\|p_0\|_{\mu_w} <\infty$ if $T=\mathbb R, 1 \le r <2$. If $\delta
_n \to0$ as $n \to\infty$, then there exists $M<\infty$ such that
%
\begin{equation}
\Pi_n\{p\in\mathcal{P}\dvtx\|p-p_0\|_r\ge M\delta_n|X_1,\ldots
,X_n\}\to0\qquad\mbox{as } n\to\infty
\end{equation}
in $P_0^{\mathbb N}$-probability.
\end{theorem}

Note that the moment condition in (\ref{pr}) is void if $r \ge2$ or if
$T=[0,1]$. If $r=1$ the rate can be taken to be $\delta_n = \varepsilon
_n$ or, more generally, $\delta_n=\gamma_n \varepsilon_n$. For $r=\infty
$ one only has at best $\delta_n = \sqrt n \varepsilon_n^2$, which is
always slower than $\varepsilon_n$ (since $\sqrt n \varepsilon_n \to
\infty$). In case $1<r<\infty$ the rate interpolates between these two
rates without, however, requiring $p_0\in L^{\infty}$.

In the case
where $p_0$ is bounded, and if it is known that the posterior
concentrates on a fixed sup-norm ball with probability approaching one,
we can refine the rates in the above theorem for $1<r<\infty$, and
retrieve the (in applications of the theorem often optimal) rate
$\varepsilon_n$ for $1 \le r \le2$.
The following theorem can be applied with $\gamma_n=1$  $\forall n$, in
which case conditions (a) and (b) require the rate $\varepsilon_n$ to
be fast enough (which in applications typically entails that a minimal
degree of smoothness of $p_0$ has to be assumed).
\begin{theorem}\label{generalbd}
Let $T, \mathcal{P}, \Pi_n$ be as in Theorem \ref{general}. Let $1< r <
\infty$, and let $\varepsilon_n \to0$ as $n \to\infty$ be a sequence
of positive numbers such that $\sqrt n \varepsilon_n \to\infty$ as $n
\to\infty$. Let $\bar r = \max(r,2)$, and set
%
\begin{equation}\label{ratebd}
\delta_n= \varepsilon_n (n \varepsilon_n^2)^{{1}/{2}-{1}/{\bar
r}} \gamma_n
\end{equation}
for some sequence $\gamma_n \ge1$. Assume either:\vspace*{8pt}

\textup{(a)} that $1<r<2$ and that $\varepsilon_n=O(\gamma_n (n \varepsilon
_n^2)^{1/r-1} )$ or

\textup{(b)} that $2 \le r < \infty$ and that $ \varepsilon_n^2 =O(\gamma
_n/\sqrt n)$.\vspace*{8pt}

Let $J_n, \mathcal P_n$ be defined as in Theorem \ref{general}, assume
that conditions (1) and~(2) in that theorem are satisfied, and that, in
addition,\vspace*{8pt}

\textup{(3)} there exists $0<B<\infty$ such that
\[
\Pi_n(p \in\mathcal P\dvtx \|p\|_\infty> B |X_1, \ldots, X_n) \to0
\]
as $n \to\infty$ in $P_0^{\mathbb N}$-probability.\vspace*{8pt}

Let $p_0 \in L^\infty(T)$ be s.t. $\|K_{J_n}(p_0)-p_0\|_r =O(\delta_n)$
and such that $\|p_0\|_{\mu_w}<\infty$ for some $w>(2-r)/r$ if
$T=\mathbb R, 1 \le r<2$. If $\delta_n \to0$ as $n \to\infty$, then
there exists $M<\infty$ s.t.
%
\begin{equation}
\Pi_n\{p\in\mathcal{P}\dvtx\|p-p_0\|_r\ge M\delta_n|X_1,\ldots
,X_n\}\to0 \qquad\mbox{as } n\to\infty
\end{equation}
in $P_0^{\mathbb N}$-probability.
\end{theorem}

\subsection{$L^r$-norm inequalities}

A main step in the proof of Theorems \ref{general} and \ref{generalbd}
[see (\ref{h}) below] is the construction of nonparametric tests for
$L^r$-alternatives, $1 \le r \le\infty$, that have sufficiently good
exponential bounds on the type-two errors. For this we first derive
sharp concentration inequalities for $L^r$-norms of centered density
estimators. It is convenient to observe that the degree of
concentration of a kernel-type density estimator around its expectation
in~$L^r$ depends on $r$, as can already be seen from comparing the
known cases $r=1, \infty$ in \cite{GN08,GG02} for kernel estimators and
\cite{GN09} for wavelets. These results are derived from Talagrand's
inequality \cite{T96} for empirical processes: let $X_1, \ldots, X_n$ be
i.i.d. with law $P$ on a measurable space $(S, \mathcal S)$, let
$\mathcal F$ be a $P$-centered (i.e., $\int f \,dP = 0$ for all $f \in
\mathcal F$) countable class of real-valued measurable functions on
$S$, uniformly bounded by the constant $U$, and set $\|H\|_\mathcal
F={\sup_{f \in\mathcal F}}|H(f)|$ for any $H\dvtx \mathcal F \to\mathbb R$.
Let $\sigma$ be any positive number such that $\sigma^2\ge\sup_{f\in
\mathcal F} E(f^2(X))$, and set $V:=n\sigma^2+2U E\|{\sum_{j=1}^n
f(X_j)}\|_\mathcal F$. Then, Bousquet's \cite{B03} version of
Talagrand's inequality, with constants, is as follows (see Theorem 7.3
in \cite{B03}): for every $x \geq0, n \in\mathbb N$,
%
\begin{equation}\label{bousq}
\Pr\Biggl\{\Biggl\| \sum_{j=1}^n f(X_j)\Biggr\|_\mathcal F \geq E \Biggl\| \sum_{j=1}^n
f(X_j)\Biggr\|_\mathcal F + \sqrt{2Vx} + Ux/3 \Biggr\} \leq2e^{-x}.
\end{equation}
This applies to our situation as follows: let $X_1, \ldots, X_n$ be
i.i.d. with density $p_0$ on $T$ with respect\vspace*{1pt} to Lebesgue measure
$\lambda$, $dP_0 = p_0 d\lambda$, and let $\hat p_n(j)=\frac{1}{n} \sum
_{i=1}^n K_j(\cdot, X_i)$ be a kernel-type estimator with $K_j$ as in
Condition \ref{kernel}. Its expectation equals $P^n_0\hat
p_n(j)(x)=EK_j(x,X) = K_j(p_0)(x)$, and we wish to derive sharp
exponential bounds for the quantity $\|\hat p_n(j)-K_j(p_0)\|
_r$ for $1 \le r \le\infty$. In case $r=\infty$ this can be achieved
by studying the empirical process indexed by
\[
\mathcal K = \{K_j(x,\cdot) - K_j(p_0)(x)\dvtx x \in T \},
\]
and in case $r<\infty$ we shall view $\hat p_n(j)-P^n_0 \hat p_n(j)$ as
a sample average of the centered $L^r(T)$-valued random variables
$K_j(\cdot, X_i) - K_j(p_0)$, and reduce the problem to an empirical
process as follows: let $s$ be conjugate to $r$, that is, $1=1/s+1/r$.
By the Hahn--Banach theorem, the separability of $L^r(T)$ implies that
there is a countable subset $B_0$ of the unit ball $B$ of
$L^s(T)$ such that
\[
\|H\|_r = \sup_{f \in B_0} \biggl|\int_\mathbb R H(t)f(t)\,dt \biggr|
\]
for all $H \in L^r(T)$. We thus have $\|\hat p_n(j) - P_0^n\hat
p_n(j) \|_r = \|P_n -P_0\|_{\mathcal K}$, where $P_n = \sum
_{i=1}^n \delta_{X_i}/n$ is the empirical measure, and where
\[
\mathcal K = \biggl\{x \mapsto\int_T f(t)K_j(t,x)\,dt - \int_T f(t)
K_j(p_0)(t)\,dt \dvtx f \in B_0 \biggr\}.
\]
To apply (\ref{bousq}) with the countable class $\mathcal K$ we need to
find suitable bounds for the envelope $U \ge{\sup_{k \in\mathcal
K}}|k(x)|$ and the weak variances $\sigma^2 \ge\sup_{k \in\mathcal
K}Ek^2(X)$. We will also apply (\ref{bousq}) in the case $r=\infty$,
and note that the corresponding empirical process suprema are over
countable subsets $B_0$ of $T$, by the continuity property of $K$ in
the convolution kernel case, and by finiteness of the $p$-variation of
the scaling function in the wavelet case (Remark~2
in~\cite{GN09}).\looseness=1

\subsubsection{Envelope and variance bounds for $\mathcal K$}

We first consider Condition~\ref{kernel}(a), the convolution kernel
case: let us write in abuse of notation $K_j(\cdot)= 2^jK(2^j\cdot)$
and $f=\delta_y, y \in B_0 \subset T$ for $r = \infty$. (One naturally
replaces~$L^s$ by the Banach space of finite signed measures if
$r=\infty$ in the arguments below.) The class $\mathcal K$ then equals
\[
\mathcal K =\bigl\{x \mapsto K_j \ast f (x) - E \bigl(K_j \ast f(X)\bigr) \dvtx f
\in B_0 \bigr\}.
\]
The bound for the envelope is seen to be of size $2^{j(1-1/r)}$: by H\"
{o}lder's inequality
%
\begin{equation} \label{env}
\|K_j \ast f\|_\infty \le \|K_j\|_r \|f\|_s \le C(K, r) 2^{j(1-1/r)}
\equiv U,
\end{equation}
a bound that remains true when $r=\infty$ since $|2^jK(2^j(x-y))| \le\|
K\|_\infty2^j$. To bound the variances, for densities $p_0 \in L^r$,
we have
%
\begin{equation} \label{var}
E(K_j \ast f)(X)^2\le\|p_0\|_r \|K_j \ast f\|_{2s}^2 \le C'(K, r) \|p_0\|_r
2^{j(1-1/r)} \equiv\sigma^2
\end{equation}
from H\"{o}lder's inequality and since $\|K_j \ast f\|_{2s}$, for $f
\in L^s$ is bounded up to constants by $2^{j(1/2-1/2r)}$, by using
Young's inequality $\|h \ast g\|_t \le\|h\|_p \|g\|_q$ for
$1+1/t=1/p+1/q, 1\le p,q,t \le\infty$.

The last estimate can be refined if $p_0$ is known to be bounded, where
we recall that $\bar r = \max(r,2)$, to
yield
%
\begin{equation} \label{varbd}
E(K_j \ast f)(X)^2\le C (p_0) 2^{j(1-2/\bar r)} \equiv\sigma^2,
\end{equation}
where $C(\cdot)$ is bounded on uniformly
bounded sets of densities. To see this, consider first $r \ge2$ and thus $s \le2$: then Young's
inequality gives, as above,
\[
E(K_j \ast f)(X)^2\le\|p_0\|_\infty\|K_j \ast f\|_{2}^2 \le C \|p_0\|
_\infty2^{j(1-2/r)} = \sigma^2.
\]
If $1<r<2$, then $p_0 \in L^\infty\cap L^1 \subset L^{s/(s-2)}$, so by
H\"{o}lder's inequality
\[
E(K_j \ast f)(X)^2 \le\|K_j \ast f\|_s^{2} \|p_0\|_{s/(s-2)} \le
C(p_0) \|K_j\|^2_1 \|f\|^2_s \le C(p_0, K).
\]
For Condition \ref{kernel}(b), so in the multiresolution case for
$T=\mathbb R$, the arguments as in (a) and obvious modifications give
the same bounds for $U, \sigma$ in view of the estimate $|{\int
_\mathbb R K_j(x,y)f(y)\,dy}| \le\Phi_j \ast|f|(x)$, which allows
us to compare wavelet projections to convolutions and proceed as above.

For Condition \ref{kernel}(c), note that, by the comments following the
statement of Condition \ref{kernel}, the\vspace*{-1pt} projection kernels have the
form $K_j=K_j^0+K_j^1+\tilde K_j$ where $\tilde K_j(x,t)=2^j\tilde
K(2^jt,2^jx)$ with $\tilde K$ majorized by a convolution kernel.
Therefore the envelope and variance bounds for the previous two cases
apply as well to this ``interior part'' of the kernel. For the boundary part,
%
\begin{equation}\label{bdary}
K_j^i(x,t)=\sum_{k=0}^{N-1}2^j\phi_k^i(2^jx)\phi_k^i(2^jt),\qquad i=0,1,
j\ge J_0,
\end{equation}
with $N$ finite and $\phi_k^i$ bounded and with bounded support, it is
immediate to check, just using H\"older's inequality, that for $f\in B_0$,
\[
\biggl\|2^j\phi_k^i(2^jx)\int_0^1\phi_k^i(2^jt)f(t)\,dt\biggr\|_\infty\le\|
\phi_k^i\|_\infty\|\phi_k^i\|_r2^{j(1-1/r)},\qquad 1\le r\le
\infty,
\]
and that
\[
2^{2j}E(\phi_k^i(2^jX))^2\biggl(\int_0^1|\phi
_k^i(2^jt)||f(t)|\,dt\biggr)^2\le\|p_0\|_r\|\phi_k^i\|_{2s}^2\|\phi_k^i\|
_r^22^{j(1-1/r)}
\]
for $p_0\in L^r$, with the refinement $\|p_0\|_\infty\|\phi_k^i\|_2^2\|
\phi_k^i\|_r^22^{j(1-2/\bar r)}$ if $\|p_0\|_\infty<\infty$. This shows
that the bounds for $U, \sigma^2$ from (a), (b) apply to (c) as well.

\subsubsection{Application of Talagrand's inequality}

To apply Talagrand's inequality we need a bound on the moment of the
supremum of the empirical process involved, provided in the following
lemma, known
for the cases $r = \infty$ (see \cite{GG02,GN09,LN10})
and, implicitly, $1 \leq r \leq 2$ (see \cite{GM07}).
As the proof is standard but somewhat lengthy it is given in the
supplementary file for this paper,~\cite{suppA}.
\begin{lemma} \label{momr}
Assume Condition \ref{kernel}\textup{(a)}, \textup{(b)} or \textup{(c)}
and that $p_0\in L^r(T)$. If $1\le r< 2$ in the cases \textup{(a)} or
\textup{(b)}, assume further that $p_0 \in L^1(\mu_s)$ for some
$s>(2-r)/r$. Then, if $1\le r<\infty$, there exists $L_r$ such that,
for all $j\ge0$ if $r\le2$, and for all $j$ such that $2^j<n$ for
$r>2$, we have
%
\begin{equation}\label{lrbounds}
E\|n(P_n-P_0)\|_{\mathcal{K}}=E\Biggl\|\sum_{i=1}^n\bigl(K_j(\cdot,
X_i)-EK_j(\cdot, X)\bigr)\Biggr\|_r\le L_r\sqrt{2^jn}.
\end{equation}
If $r=\infty$, for $p_0$ and $\Phi$ bounded, there exists a
constant $L_\infty$ such that for all~$j$ satisfying $2^jj<n$ we have
%
\begin{equation}\label{linfbound}\quad
E\|n(P_n-P_0)\|_{\mathcal{K}}=E\Biggl\|\sum_{i=1}^n\bigl(K_j(\cdot,
X_i)-EK_j(\cdot, X)\bigr)\Biggr\|_\infty\le L_\infty\sqrt{2^jjn}.
\end{equation}
\end{lemma}

We are now ready to apply (\ref{bousq}): for $V= n\sigma^2 + 2U E
\|\hat p_n(j)-E \hat p_n(j) \|_r$ we have the bound
\begin{eqnarray*}
\Pr\biggl\{n\|\hat p_n(j)-P_0^n \hat p_n(j) \|_r &\ge&
nE\|\hat p_n(j)-P_0^n \hat p_n(j)
\|_r
+ \sqrt{2Vx} + \frac{Ux}{3} \biggr\} \le2e^{-x}.
\end{eqnarray*}
This can be further simplified, using the standard inequalities $\sqrt
{a+b} \le\sqrt a + \sqrt b, \sqrt{ab} \le(a+b)/2$, to
\begin{eqnarray*}
&&\Pr\bigl\{n\|\hat p_n(j)-P_0^n\hat p_n(j) \|_r \ge \tfrac
{3}{2} nE\|\hat p_n(j)-P_0^n \hat p_n(j) \|_r
+ \sqrt{2n \sigma^2 x} + \tfrac{7}{3}Ux \bigr\} \\
&&\qquad\le2e^{-x}.
\end{eqnarray*}
Combining the moment estimate Lemma \ref{momr} with (\ref{env}) and
(\ref{var}), we obtain, for $2^j j(r)<n$ with $j(\infty)=j$ and
$j(r)=1$ for $r<\infty$,
%
\begin{eqnarray} \label{fintal}
&&\Pr\bigl\{n\|\hat p_n(j)-P_0^n\hat p_n(j) \|_r \nonumber\\[-8pt]\\[-8pt]
&&\qquad\ge C
\bigl(\sqrt{2^jn j(r)}
+ \sqrt{n 2^{j(1-1/r)} \|p_0\|_r x} + 2^{j(1-1/r)}x \bigr)
\bigr\}
\le2e^{-x}\nonumber
\end{eqnarray}
for some constant $C$, and in the case where $\|p_0\|_\infty<\infty$ we
have, analogously, from (\ref{varbd}),
%
\begin{eqnarray} \label{fintalbd}
&&\Pr\bigl\{n\|\hat p_n(j)-P_0^n\hat p_n(j) \|_r \nonumber\\[-8pt]\\[-8pt]
&&\qquad\ge C
\bigl(\sqrt{2^jn j(r)}
+ \sqrt{n 2^{j(1-2/\bar r)} \|p_0\|_\infty x} + 2^{j(1-1/r)}x \bigr)
\bigr\} \le2e^{-x}.\nonumber
\end{eqnarray}
If we take $\varepsilon_n$, $\delta_n$,
$2^{j_n} \sim n \varepsilon_n^2$ as in Theorems \ref{general}, \ref
{generalbd}, and if $\|p_0\|_r$ is bounded by a fixed constant $B$,
then the choice $x = Ln \varepsilon_n^2$ gives for every $L$ and
$M=M(L,K,B)$ large enough, after some simple computations using the
conditions on $\varepsilon_n, \delta_n$ from the theorem, that
\[
nM\delta_n \ge C\bigl(\sqrt{2^{j_n} j_n(r)n} + \sqrt{\|p_0\|_r n
2^{j_n(1-1/r)} L n \varepsilon_n^2} + 2^{j_n(1-1/r)}Ln \varepsilon_n^2
\bigr)
\]
and, likewise, if $\|p_0\|_\infty$ is bounded by a fixed constant, the
corresponding choice of $\delta_n, M$ also satisfies
\[
nM\delta_n \ge C\bigl(\sqrt{2^{j_n} j_n(r)n} + \sqrt{C(p_0) n
2^{j_n(1-2/\bar r)} L n \varepsilon_n^2} + 2^{j_n(1-1/r)} Ln \varepsilon
_n^2 \bigr).
\]
Moreover for $\|p_0\|_r \ge\zeta>0$ we have
\[
n \|p_0\|_r \ge C\bigl(\sqrt{2^{j_n} j_n(r)n} + \sqrt{\|p_0\|_r n
2^{j_n(1-1/r)} L n \varepsilon_n^2} + 2^{j_n(1-1/r)}Ln \varepsilon_n^2
\bigr)
\]
from some index $n_0$ onwards that depends only on $C, \zeta$.

Using these inequalities in (\ref{fintal}), (\ref{fintalbd}), we
conclude that in both cases, for every $0<L<\infty$ we can find a large
enough $M(L,K,B)$ such that
%
\begin{equation} \label{taluse}
\Pr\{n\|\hat p_n(j_n)-P_0^n\hat p_n(j_n) \|_r \ge M
n\delta_n \} \le2 e^{-Ln \varepsilon_n^2}
\end{equation}
and, likewise, for $n$ large enough,
%
\begin{equation} \label{talusep}
\Pr\{n\|\hat p_n(j_n)-P_0^n \hat p_n(j_n) \|_r \ge
n\|p_0\|_r/3 \} \le2e^{-Ln \varepsilon_n^2}.
\end{equation}

\subsection{\texorpdfstring{Proof of Theorems \protect\ref{general} and \protect\ref{generalbd}}
{Proof of Theorems 2 and 3}}

Using the small ball estimate from condition (2), it suffices to
construct tests (indicator functions) $\phi_n=\phi_n(X_1,\ldots,\allowbreak X_n;p_0)$ such that
%
\begin{eqnarray}\label{h}
&\displaystyle P_0^n\phi_n\to0 \qquad\mbox{as } n\to\infty  \quad\mbox{and}&\nonumber\\[-8pt]\\[-8pt]
&\displaystyle \sup_{p\in\mathcal
{P}_n\dvtx\|p-p_0\|_r \ge
M\delta_n}P^n(1-\phi_n)\le2e^{-(C+4)n\varepsilon_n^2}&\nonumber
\end{eqnarray}
for $n$ large enough; see the proof of Theorem 2.1 in \cite{GGV00}.

Consider first Theorem \ref{general}. Let $\hat p_n$ be a kernel-type
density estimator based on an i.i.d. sample $X_1,\ldots, X_n$ of common
law $P_0$, $n\in\mathbb N$, at resolution $J_n$. For $M_0$, a~constant
to be chosen below, set
$T_n=\|\hat p_n-p_0\|_r$ and $\phi_n=I(T_n>M_0 \delta_n)$. Note
that $\phi_n$ is the (indicator of the) rejection region of a natural
test of the hypothesis $H_0\dvtx p=p_0$. Then we have
\begin{eqnarray*}P_0^n\phi_n&=&P_0^n\{\|\hat p_n-p_0\|_r>M_0\delta
_n\}\\
&\le& P_0^n\{\|\hat p_n-P_0^n\hat p_n\|_r>M_0\delta_n-\|P_0^n\hat
p_n-p_0\|_r \}.
\end{eqnarray*}
Since $\|K_{J_n}(p_0)-p_0\|_r \le c'\delta_n$ for some $c'>0$ by
assumption, we have for all $n$ large enough,
$P_0^n\phi_n\le P_0^n\{\|\hat p_n-P_0^n\hat p_n\|_r>(M_0-c')\delta
_n\}$. Then using inequality~(\ref{taluse}), we have for some
constant $L_1$
for some constant $L_1$, choosing~$M_0$ large enough, that, as $n \to\infty$,
%
\begin{equation}\label{h0}
P_0^n\phi_n\le2e^{-L_1n \varepsilon_n^2} \to0.
\end{equation}

Let now $p$ be a density in $\mathcal{P}_n$ such that $\|p-p_0\|_r \ge
M\delta_n$ (the alternatives). Set $dP(x)=p(x)\,dx$. We have, from the
triangle inequality,
%
\begin{eqnarray}\label{h11}\qquad
P^n(1-\phi_n) &=& P^n\{\|\hat p_n-p_0\|_r \le M_0\delta_n\}
\nonumber\\
&\le& P^n\{\|\hat p_n-P^n\hat p_n\|_r \ge\|p-p_0\|_r-M_0\delta
_n-\|P^n\hat p_n-p\|_r\} \\
&\le& P^n\bigl\{\|\hat p_n-P^n\hat p_n\|_r \ge\|p-p_0\|
_r-\bigl(M_0+C(K)\bigr)\delta_n \bigr\}\nonumber
\end{eqnarray}
since by assumption on $\mathcal{P}_n$, $\sup_{p \in\mathcal P_n}\|
P^n\hat p_n-p\|_r\le C(K)\delta_n$,
uniformly in \mbox{$p\in\mathcal{P}_n$}.

To complete the estimation of the last probability, we consider first
$r>1$. For those $p \in\mathcal P_n$ satisfying $\|p\|_r\ge2\|p_0\|
_r$ we have $\|p-p_0\|_r \ge\|p\|_r/2 \ge\|p_0\|_r$, and, using
inequality (\ref{talusep}) for $p_0=p$, we deduce, that for all $L>0$,
there exists $n_0\in\mathbb N$ such that for all $n\ge n_0$,
%
\begin{eqnarray}\label{h12}
&&\sup_{p\in\mathcal{P}_n\dvtx \|p\|_r\ge2\|p_0\|_r}P^n(1-\phi_n)\nonumber\\
&&\qquad\le \sup
_{p\in\mathcal{P}_n, \|p\|_r\ge2\|p_0\|_r}P^n\biggl\{\|\hat p_n-P^n\hat
p_n\|_r>\frac{\|p\|_r}{3}\biggr\}\\
&&\qquad \le 2e^{-Ln\varepsilon_n^2}.\nonumber
\end{eqnarray}
For those $p \in\mathcal P_n$ for which $\|p\|_r < 2\|p_0\|_r$, we
apply (\ref{taluse}) with $p=p_0$ and use as well $\|p-p_0\|_r\ge M
\delta_n$ to obtain that for all $L>0$ there exists $M$ large enough
such that
%
\begin{eqnarray}\label{h13}
&&\sup_{p\in\mathcal{P}_n\dvtx \|p\|_r< 2\|p_0\|_r, \|p-p_0\|_r \ge M\delta
_n}P^n(1-\phi_n) \nonumber\\
&&\qquad \le\sup_{p\in\mathcal{P}_n\dvtx \|p\|_r< 2\|p_0\|_r, \|p-p_0\|_r\ge
M\delta_n}P^n\bigl\{\|\hat p_n-P^n\hat p_n\|_r\nonumber\\[-8pt]\\[-8pt]
&&\qquad\hspace*{164.8pt}>\bigl(M-M_0-C(K)\bigr)\delta
_n\bigr\}\nonumber\\
&&\qquad \le2e^{-Ln\varepsilon_n^2}.\nonumber
\end{eqnarray}
We conclude from (\ref{h11}) and (\ref{h12}) that for any $L>0$ there
exists $n_L<\infty$ such that
%
\begin{equation}\label{h14}
\sup_{p\in\mathcal{P}_n\dvtx \|p-p_0\|_r\ge M\delta_n}P^n(1-\phi_n)\le
2e^{-Ln\varepsilon_n^2}.
\end{equation}
Now (\ref{h0}) and (\ref{h14}) prove (\ref{h}) if $r>1$. If $r=1$ the
above case distinction is not necessary as $\|p\|_1=1$ always holds, so
that the proof of the second case applies with the full supremum over
$\{p \in\mathcal P_n\dvtx \|p-p_0\|_1 \ge M\delta_n\}$. This completes the
proof of Theorem \ref{general}.

To prove Theorem \ref{generalbd} we argue similarly, and only have to
slightly modify the derivation of the error probabilities of the tests:
when it is known that the posterior concentrates on a fixed sup-norm
ball of radius $B$, then we can restrict the alternatives in (\ref{h})
further to densities bounded by $B$, and, using (\ref{taluse}) with
$p=p_0$ and the present choice of $\delta_n$, we also obtain
\begin{eqnarray*}
&&\sup_{p\in\mathcal{P}_n\dvtx \|p\|_\infty\le B, \|p-p_0\|_r \ge M\delta
_n}P^n(1-\phi_n) \\
&&\qquad \le\sup_{p\in\mathcal{P}_n\dvtx \|p\|_\infty\le B, \|p-p_0\|_r\ge
M\delta_n}P^n\bigl\{\|\hat p_n-P^n\hat p_n\|_r>\bigl(M-M_0-C(K)\bigr)\delta
_n\bigr\}\\
&&\qquad \le2e^{-Ln\varepsilon_n^2}.
\end{eqnarray*}

\section{Remaining proofs}

\subsection{\texorpdfstring{Proofs of Propositions \protect\ref{unif}, \protect\ref{dir} and \protect\ref{dir1}}
{Proofs of Propositions 1, 2 and 3}}

\mbox{}

\begin{pf*}{Proof of Proposition \ref{unif}}
Since $\|U_{\alpha}\|_\infty\le C$ almost surely for some fixed
constant $C=C(B,\alpha, \psi)$, we infer $\|p^{U, \alpha}\|_{\alpha, r,
\infty} \le D(B, \alpha, \psi)$ almost surely for $1 \le r \le\infty$.
In particular the prior is supported in a ball of bounded densities,
hence so is the posterior, and we can attempt to apply Theorems~\ref{general}
(for $r=1,\infty$) and \ref{generalbd} for
($1<r<\infty$), which we shall do with the choice $\varepsilon_n =
(n/\log n)^{-\alpha/(2\alpha+1)}$.

We verify the small ball estimate in the second condition in
Theorem~\ref{general}. By Lemma 3.1 in \cite{VV08} we can lower bound the
prior\vadjust{\goodbreak}
probability in question by $\Pr\{\|{\log p_0-U_\alpha}\|_\infty\le
c\varepsilon_n \}$ for some constant $c>0$. Since
\[
\|h\|_\infty\le C(\phi,\psi) \max\biggl(\sup_k |\alpha_k(h)|, \sum
_\ell\sup_k 2^{\ell/2}|\beta_{\ell k}(h)| \biggr)
\]
for any continuous function $h$ on $[0,1]$ and some constant $C(\phi,
\psi)$, we can lower bound the last probability, writing $\alpha_k,
\beta_{\ell k}$ for the wavelet coefficients of $\log p_0$, by
\begin{eqnarray*}
\hspace*{-7pt}&&
\Pr\biggl\{\max\biggl(\sup_{k=0, \ldots, N} |\alpha_k-u_{0k}|, \sum_\ell
\sup_k 2^{\ell/2}\bigl|\beta_{\ell k}-2^{-\ell(\alpha+1/2)}u_{\ell k}\bigr|
\biggr) \le c'\varepsilon_n \biggr\}
\\
\hspace*{-7pt}&&\quad=
\Pr\Bigl\{ {\max_k} |\alpha_k-u_{0k}| \le c' \varepsilon_n\Bigr\} \Pr
\biggl\{\sum_{\ell\ge J_0} \max_{k\le2^\ell} 2^{\ell/2}\bigl|\beta_{\ell k}
-2^{-\ell(\alpha+1/2)}u_{\ell k}\bigr| \le c'\varepsilon_n \biggr\},
\end{eqnarray*}
where $N, J_0$ depend only on the wavelet basis (see before Definition
\ref{besov}). Since $|\alpha_k| \le B$ and since the $u_{0k}$ are
$U(-B,B)$, the first\vspace*{1pt} probability exceeds $(c' \varepsilon_n/\allowbreak2B)^{N+1}=
e^{-(N+1)\log(2B/c'\varepsilon_n)}$ which is bounded below by
$e^{-c\log(1/\varepsilon_n)}$ for some $c>0$ that depends only on $B$,
$\alpha$ and the wavelet basis. For the second probability set $b_{\ell
k}\equiv2^{\ell(\alpha+1/2)}\beta_{\ell k} , \ell\ge J_0$, and
$M(J)\equiv\sum_{\ell=J_0}^J \sum_ {k=0}^{2^\ell-1} 1 \le2\cdot2^J$,
and note that $|b_{\ell k}|\le\|{\log p_0}\|_{\alpha,\infty}\le B$.
Choosing $J=J_n\ge J_0$ large enough and of order $\varepsilon_n \simeq
2^{-J\alpha}$, this probability is bounded below by
\begin{eqnarray*}
&&\Pr\Biggl\{{\sum_{\ell=J_0}^J 2^{-\ell\alpha} \sup_k} |b_{\ell
k}-u_{\ell k}| \le c'\varepsilon_n - C(\psi,B)2^{-J\alpha} \Biggr\}\\
&&\qquad \ge\Pr\Bigl\{{\max_{\ell\le J}\max_{k\le2^\ell}}|b_{\ell k}-u_{\ell
k}| \le c'' \varepsilon_n \Bigr\} \\
&&\qquad= \prod_{\ell\le J}\prod_{k\le2^\ell} \Pr\{|b_{\ell
k}-u_{\ell k}| \le c'' \varepsilon_n\} \ge\biggl(\frac
{c''\varepsilon_n}{2B}\biggr)^{M(J)} \\
&&\qquad\ge e^{-c'''{\log
(1/\varepsilon_n)}/{\varepsilon^{1/\alpha}_n}}
\end{eqnarray*}
for $n$ large enough and some $c'''>0$ that depends only on $B$, $
\alpha$ and the wavelet basis. Summarizing we have, by definition of
$\varepsilon_n$, that the $\Pi^\alpha$ probability in condition (2) of
Theorem \ref{general} is bounded from below by
%
\begin{equation} \label{lowp}\hspace*{4pt}
\Pr\{\|{\log p_0-U_\alpha}\|_\infty\le c\varepsilon_n \}\ge
e^{-c\log(1/\varepsilon_n)} e^{-c'''\log(1/\varepsilon_n)/\varepsilon
^{1/\alpha}_n} \ge e^{-Cn \varepsilon_n^2}
\end{equation}
for some $C$ that depends only on $B$, $\alpha$ and the wavelet basis,
which proves that condition (2) holds.

We next verify the bias condition with $\mathcal P_n = \operatorname{supp} (\Pi)$ so
that $\Pi(\mathcal P \setminus\mathcal P_n)=0$. We bound the
$L^r$-norm of the approximation errors of any element in~$\mathcal P_n$
by a constant times $\delta_n$, where we take $\gamma_n$ equal to $\log
n$ to a sufficiently large\vadjust{\goodbreak} power chosen below. Since $2^{J_n}\ge c n
\varepsilon_n^2 \ge c n^{1/(2\alpha+1)} $ we have, using
$B_{r1}^0([0,1]) \subset L^r([0,1])$ and $p \in\mathcal C^\alpha([0,1])$,
\[
\|K_{J_n}(p)-p\|_{r} \le c \sum_{\ell=J_n}^\infty2^{\ell({1}/{2}-
{1}/{r})} \Biggl(\sum_{k=1}^{2^j} |\beta_{\ell k}(p)|^{r}\Biggr)^{1/r} \le
c'(B,r) \sum_{\ell=J_n}^\infty2^{-\ell\alpha},
\]
which is $O(\varepsilon_n )$, so the bias condition is satisfied for
some $C(K)$ large enough, both for $\mathcal P_n$, as well as for $p_0$.

Finally condition (c) from Theorem \ref{general} and (a), (b) from
Theorem \ref{generalbd}, as well as $\delta_n \to0$, are verified for
this choice of $\varepsilon_n$ and under the conditions on~$\alpha,r$,
except for the cases $\alpha=0$ or $\alpha=1/2, r =\infty$, where the
result trivially follows from $\delta_n$ being bounded from below by a
constant multiple of $\log n$ (and as the prior is supported in a
$L^r$-bounded set).
\end{pf*}
\begin{pf*}{Proof of Proposition \ref{dir}}
We apply Theorem \ref{general} with $r=\infty$. We have from the proof
of Theorem 5.1 in \cite{GV01} that for $\varepsilon_n = (\log n)^\kappa
/ \sqrt n, \kappa\ge1$, the small-ball estimate in condition (2) of
Theorem \ref{general} is satisfied. Choose $\gamma_n$ in such a way
that $\delta_n$ equals $(\log n)^\eta/\sqrt n$ where $\eta>\kappa$.
For the bias, we take $\mathcal P_n$ to be the support of $\Pi$ and
consider a Meyer-wavelet basis and the wavelet projection onto it, with
$2^{J_{n}} = c(\log n)^{2 \kappa}$, where $c$ is a large enough
constant that depends on $\inf\{\sigma\dvtx \sigma\in \operatorname{supp}(G)\}$,
and apply Proposition 4 in \cite{LN10} with $s=2$ and suitable $\tilde
c_0$, to see that $\|K_{J_n}(p_{F, \sigma})-p_{F, \sigma}\|_\infty=
o(1/n)$ uniformly in the support of $\Pi$. A more detailed proof is in
the supplementary file \cite{suppA}.\looseness=-1
\end{pf*}

\begin{pf*}{Proof of Proposition \ref{dir1}}
Taking $\varepsilon_n = M'(n/\log n)^{-\alpha/(2\alpha+1)}$, and noting
$\varepsilon_n^{-1/\alpha} =O(n \varepsilon_n^2)$, we can take $J_n$
such that $2^{j_n} \le2^{J_n} \le c n \varepsilon_n^2$ for every~$n$,
some \mbox{$c>0$}. Taking $K(x,y)$ equal to the Haar wavelet projection kernel
(CDV-wavelet of regularity $S=0$), we conclude that $\|K_{J_n}(p)-p\|
_r=0$ $\Pi_{j_n}$-a.s. $\forall n$, so condition (1) in Theorem \ref
{general} is satisfied with $\mathcal P_n$ equal to the support of $\Pi
_{j_n}$. The small ball estimate (2) follows, as in the proof of
Theorem 1 (\cite{S07}, pages 636 and 637, with $k_0 = 2^{j_n}$, and
approximating $p_0$ by $K_{j_n}(p_0)$ s.t. $\|K_{j_n}(p_0)-p_0\|_1 \le
\varepsilon_n/2$ for $M'$ large enough), and from the second inequality
in (\ref{lowp}). The bias condition for $p_0$ is satisfied by standard
approximation properties of Haar wavelets. The result now follows from
first applying Theorem \ref{general} with $r=1, \infty$ and then using
the conclusion that the posterior concentrates on a \mbox{$\|\cdot\|_\infty$}
neighborhood of $p_0$ to invoke Theorem \ref{generalbd} for the cases
$1<r<\infty$.
\end{pf*}

\subsection{\texorpdfstring{Proof of Proposition \protect\ref{fbmp}}{Proof of Proposition 4}}

We shall construct subsets of $\mathcal P$ on which we can control the
approximation errors from (\ref{pr}). We define H\"{o}lder spa\-ces.~For~$\alpha, \tau\ge0$ positive\vadjust{\goodbreak} real numbers, define the norm
$\|f\|_{\alpha,\infty,\tau}:=\break{\sum_{k=0}^{[\alpha]}}\|f^{(k)}\|_\infty+
H(\alpha, \tau, f)$ where
\[
H(\alpha, \tau, f)=\sup_{0<t<1}\frac{{\sup_{h\dvtx|h|\le t, x+h \in
[0,1]}\sup_{x\in[0,1]}}|f^{(k)}(x+h)-f^{(k)}(x)|}{t^{\{\alpha\}}(\log
t^{-1})^{\tau}},\vadjust{\goodbreak}
\]
and where we take $\|f^{(k)}\|_\infty=\infty$ if $f^{(k)}$ does not
exist. Define, moreover, $C^{\alpha,\infty,\tau}([0,1]):=\{f\dvtx
[0,1]\to\mathbb R\dvtx \|f\|_{\alpha,\infty,\tau}<\infty\}$. The case
$\tau=0$ specialises to the strict $\alpha$-H\"older case $C^\alpha([0,1])$.

In case $1\le r <\infty$, we shall use approximation theoretic
properties of the reproducing kernel Hilbert spaces (RKHSs) of
$B_\alpha, \bar B_\alpha$, which are Sobolev spaces. Recall that the
RKHS $\mathbb H(1/2)$ of Brownian motion on $[0,1]$ is the space of
absolutely continuous functions that are zero at zero and whose first
derivatives are in $L^2([0,1])$, equipped with the inner product
$\langle f,g\rangle_{\mathbb H(1/2)}=\int_0^1f'g'$. Then, the RKHS of
integrated Brownian motion $B_\alpha$ is
\[
\mathbb
H(\alpha)=\biggl\{\int_0^t\int_0^{[\alpha]-1}\cdots\int
_0^{t_1}f(s)\,ds\,dt_1\cdots
dt_{[\alpha]-1}\dvtx f\in\mathbb H(1/2)\biggr\}
\]
with inner product
$\langle f,g\rangle_{\mathbb
H(\alpha)}=\int_0^1f^{([\alpha]+1)}g^{([\alpha]+1)}$. Finally, $f\in
\bar{\mathbb H}(\alpha)$, the RKHS of $\bar B_\alpha$, iff $f=
P_{[\alpha]}+g$ where $P_\alpha$ is a polynomial of degree $[\alpha]$
and $g\in\mathbb H(\alpha)$, and note that
$P_{[\alpha]}(t)=\sum_{i=0}^{[\alpha]}f^{(i)}(0)t^i/i!$; the inner
product in $\bar{\mathbb H}(\alpha)$ is $\langle
f,g\rangle_{\bar{\mathbb
H}(\alpha)}=\sum_{i=0}^{[\alpha]}f^{(i)}(0)g^{(i)}(0)+\int
_0^1f^{([\alpha]+1)}g^{([\alpha]+1)}$;
see, for example, \cite{VV08b}. The spaces $\bar{\mathbb H}(\alpha)$
are precisely the Sobolev spaces $H^{\alpha+1/2}$, and other equivalent
norms may be used below.

We will also require the following definition. For a $B$-valued
Gaussian random vector $W$, $B$ a Banach space, and for $w\in B$, the
``concentration function'' $\phi_w^W(\varepsilon)$ of $W$ at $w$ is defined
as
%
\begin{equation}\label{conc}
e^{-\phi_w^W(\varepsilon)}=\Pr\{\|W-w\|<\varepsilon\}.
\end{equation}

The following result is a consequence of Borell's isoperimetric
inequality~\cite{B75}, and is essentially contained in the proof of
Theorem 2.1 in \cite{VV08}.
\begin{proposition}\label{borelll0}
Let $\alpha\in\{n-1/2\dvtx n \in\mathbb N\}$, denote by $\bar{\mathbb
H}_1(\alpha)$ the unit ball of $\bar{\mathbb H}(\alpha)$ and let $B^{1}
= \{f \in C([0,1])\dvtx \|f\|_{ \infty} \le1\}$. Let $\varepsilon_n$
satisfy $\phi_0^{\bar B_\alpha}(\varepsilon_n)\le n\varepsilon_n^2$ for
all $n$. Then the released integrated Brownian motion process $\bar
B_\alpha$ has a~ver\-sion, that we continue denoting by $\bar B_\alpha$,
such that for every $C>0$, $D>0$,
\[
\Pr\{\bar B_\alpha\notin M_n \bar{\mathbb H}_1(\alpha)+
\varepsilon_nB^{1} \} \le De^{-(C+4) n \varepsilon_n ^2},
\]
where $ M_n=M_n(C,D)=-2\Phi^{-1}(De^{-(C+4)n\varepsilon_n^2}
)\simeq\sqrt n \varepsilon_n$ and $\Phi$ is the standard normal
distribution function.
\end{proposition}
\begin{pf}
Borell's inequality (e.g., Theorem 4.3.3 in \cite{B98}) implies
%
\begin{equation}\label{Brkhs}
\Pr\{\bar B_\alpha\notin M_n \bar{\mathbb H}_1(\alpha)+ \varepsilon
_nB^{1} \} \le1 - \Phi(a_n+M_n),
\end{equation}
where $a_n$ solves the equation $\Phi(a_n) = \Pr\{\|\bar B_\alpha\|
_\infty\le\varepsilon_n\}\ge e^{-n\varepsilon_n^2}$. It then follows
($C+4>1$) that $a_n\ge-M_n/2$, which implies
\[
1 - \Phi(a_n+M_n)\le\Phi(-M_n/2)=D
e^{-(C+4)n\varepsilon_n^2}.\qquad\qed
\]
\noqed\end{pf}

In particular, taking $D=\Pr\{\|\bar B_\alpha\|_\infty\le c\}
$ for any $c>0$, this proposition gives
%
\begin{equation}\label{borell0}
\Pr\{\bar B_\alpha\notin M_n \bar{\mathbb H}_1(\alpha)+
\varepsilon_nB^{1} | \|\bar B_\alpha\|_\infty\le c\} \le
e^{-(C+4) n \varepsilon_n ^2}
\end{equation}
with $M_n$ depending on $C$ and $c$, and of the order $\sqrt n
\varepsilon_n$.

In case $r=\infty$ we need a different result that reflects the almost
sure H\"{o}lder regularity of the trajectories of $B_\alpha$.
\begin{proposition}\label{borell}
For all $\alpha\in\{n-1/2\dvtx n \in\mathbb N\}$, integrated Brownian
motion has a version, that we continue denoting by $B_\alpha$, with
almost all its sample paths in $C^{\alpha,\infty,1/2}([0,1])$ and for
every $D>0$ there exist $t_\alpha<\infty$ and $L_\alpha<\infty$ such that
%
\begin{equation}\label{Borell}
\Pr\{\|B_\alpha\|_{\alpha,\infty,1/2}\ge t\}\le
De^{-L_\alpha t^2},\qquad t\ge t_\alpha.
\end{equation}
The same is true for the processes $\bar B_\alpha=\sum_{k=0}^{[\alpha
]+1}Z_kt^k/k!+B_\alpha$, that is,
%
\begin{equation}\label{Borelll}
\Pr\{\|\bar B_\alpha\|_{\alpha,\infty,1/2}\ge t\}\le
De^{-L_\alpha t^2},\qquad t\ge t_\alpha,
\end{equation}
for possibly different $L_\alpha(D)$ and $t_\alpha(D)$, for all $D>0$.
\end{proposition}
\begin{pf}
By a classical\vspace*{1pt} result of L\'evy (see also Theorem IV.5 in \cite{CKR93})
Brownian motion $B_{1/2}$ has a version in $C^{1/2,\infty,1/2}([0,1])$.
Since, for $\alpha>1$, by the definitions,
\[
\|B_\alpha\|_{\alpha,\infty,1/2}=\|B_\alpha\|_\infty+\|B_\alpha'\|
_{\alpha-1,\infty,1/2}=\|B_\alpha\|_\infty+\|B_{\alpha-1}\|_{\alpha
-1,\infty,1/2},
\]
and $\|B_\alpha\|_\infty<\infty$ a.s., induction extends the result to
all $\alpha\in\{n-1/2\dvtx n \in\mathbb N\}$.

For $0<\alpha<1$, Theorem III.6 in \cite{CKR93} shows that the norms $\|
f\|_{\alpha,\infty,1/2}$ and $\|f\|_{\alpha,\infty,1/2}^{(d)}$ are
equivalent, where $\|f\|_{\alpha,\infty,1/2}^{(d)}$ is defined as
%
\begin{eqnarray}\label{norm2}
\|f\|_{\alpha,\infty,1/2}^{(d)}:\!&=&\|(y_i^f,y_{j,k}^f)\|_{\alpha,\infty
,\infty}\nonumber\\[-8pt]\\[-8pt]
&=&\sup\biggl\{|y_0^f|,|y_1^f|,\max_{k,j}\frac{2^{\alpha j}}{\sqrt
{j\log2}}|y_{j,k}^f|\biggr\}\nonumber
\end{eqnarray}
with
\begin{eqnarray}\label{y}
y_0^f&=&f(0),\nonumber\\
y_1^f&=&3^{-1/2}\bigl(f(1)-f(0)\bigr),
\\
%
y_{j,k}^f&=&(3\cdot2^J)^{-1/2}\biggl[f\biggl(\frac{2k-1}{2^{j+1}}
\biggr)-\frac{1}{2}\biggl(f\biggl(\frac{k}{2^j}\biggr)+f\biggl(\frac{k-1}{2^j}
\biggr)\biggr)\biggr]\nonumber
\end{eqnarray}
for $k=1,\ldots,2^j, j=0,1,\ldots.$ Obviously,
$\|\cdot\|_{\alpha,\infty,1/2}^{(d)}$ is a supremum norm on a sequence
space; more specifically, it is the sup of the absolute values of a~%
countable number of linear functionals on the space
$C^{\alpha,\infty,1/2}([0,1])$ (linear combinations\vadjust{\goodbreak} of point
evaluations). Hence Lemma 3.1 and inequality (3.2) in~\cite{LT91} (this
last inequality even with $\pi^2/2$ replaced by 2) apply to
$\|B_\alpha\|_{\alpha,\infty,1/2}$, giving (\ref{Borell}) for $D=1$.
For $D<1$, take $t_\alpha'\ge t_\alpha$ such that $D\ge
e^{-(L_\alpha/2)(t'_\alpha)^2}$ and $L'_\alpha=L_\alpha/2$. If
$\alpha>1$, then the result follows by applying these inequalities to
the $C^{\{\alpha\},\infty,1/2}$-norm of the $[\alpha]$th derivative of
the process and to the sup norms of the process and of its derivatives
of order smaller than $[\alpha]$. Since~(\ref{Borell}) is obviously
true for the processes $Z_kt^k$, it is true as well for $\bar B_\alpha$
possibly with a different constant, which gives (\ref{Borelll}).
\end{pf}

Again, taking $D=\Pr\{\|\bar B_\alpha\|_\infty\le c\}$, for
any $c>0$, this proposition gives
%
\begin{equation}\label{Borellll}
\Pr\{\|\bar B_\alpha\|_{\alpha,\infty,1/2}\ge t | \|\bar
B_\alpha\|_\infty\le c\}\le e^{-L_\alpha t^2},\qquad t\ge t_\alpha,
\end{equation}
$L_\alpha$ and $t_\alpha$ depending on $c$.\vadjust{\goodbreak}

These two consequences of Borell's inequality imply that the integrated
Brownian motions concentrate on suitable subsets of $C([0,1])$, and the
following lemma achieves the same for the normalized trajectories of
the processes $e^{\bar B_\alpha(t,\omega)}$.
\begin{lemma} \label{BIAS} Let $\alpha\in\{n-1/2\dvtx n \in\mathbb N\}$,
and let $K_j$ be a CDV-projection kernel of regularity $\alpha+1/2$, at
resolution $j\ge0$.\vspace*{8pt}

(1) (Case $1 \le r < \infty$.) Let $f \in \{M_n \bar{\mathbb
H}_1(\alpha)+\varepsilon_n B^{1}, \|f\|_\infty\le c\}$, where
$\bar{\mathbb H}_1(\alpha)$ is the unit ball of the RKHS of $\bar
B_\alpha$ and set $p=e^{f}/\int_0^1 e^{f}$. Then, for $\bar r = \max
(2,r)$ and some $C>0$,
\[
\|K_j(p)-p\|_r \le C\bigl( M_n 2^{-j(\alpha+1/\bar r)} + \varepsilon
_n\bigr).
\]

(2) (Case $r = \infty$.) Let $f$ satisfy $\|f\|_\infty\le c$ and $\|f\|
_{\alpha,\infty,1/2} \le L\sqrt n \varepsilon_n$, and let $p$ be as
above. Then, for some $C>0$,
\[
\|K_j(p)-p\|_\infty\le C\sqrt n \varepsilon_n 2^{-j\alpha} \sqrt j.
\]
\end{lemma}
\begin{pf}
We first consider $1 \le r < \infty$. Since $\|f\|_\infty\le c$ we
have $e^{-c}\le\int_0^1e^{f}\le e^c$ so, $\int K_j(x,y)(\cdot)(y)\,dy$
being a linear operator, it suffices to bound $\|K_j(e^{f})-e^{f}\|_r$.
Writing $f=f_1+f_2$ with $f_1 \in M_n\bar{\mathbb H}_1(\alpha)$ and $f_2
\in \varepsilon_n B^1$, we see that $\|f_2\|_\infty\le\varepsilon_n<c$, $\|f_1\|
_\infty\le c+\varepsilon_n<2c$, and in particular,
$|e^{f_2(x)}-e^{f_2(y)}|\le e^c|f_2(x)-f_2(y)|$. Note also that, for
some constant $C(K)<\infty$,
$\|2^{-j}K_j(x, x+2^{-j}\cdot)\|_1 \le C(K)$. Then we have
\begin{eqnarray*}
&&|K_j(e^{f})-e^{f}|(x) \\
&&\qquad= \biggl|\int2^{-j}K_j(x, x+2^{-j}u)\bigl(e^{(f_1+f_2)(x+2^{-j}u)}
- e^{(f_1+f_2)(x)} \bigr) \,du \biggr| \\
&&\qquad\le \biggl|e^{f_2(x)}\int2^{-j}K_j(x,x+2^{-j}u)
\bigl(e^{f_1(x+2^{-j}u)} - e^{f_1(x)} \bigr) \,du \biggr| \\
&&\qquad\quad{} + \biggl|\int2^{-j}K_j(x,x+2^{-j}u) e^{f_1(x+u2^{-j})}
\bigl(e^{f_2(x+2^{-j}u)} - e^{f_2(x)} \bigr) \,du \biggr| \\
&&\qquad \le e^c |K_j(e^{f_1})(x)-e^{f_1}(x)| + {2e^{3c} \sup_x} \|
2^{-j}K_j(x, x+2^{-j}\cdot)\|_1 \varepsilon_n.
\end{eqnarray*}
The $L^r([0,1])$-norm of the second term is bounded by a fixed constant
times~$\varepsilon_n$, and it remains to control the $L^r([0,1])$-norm
of the first term in the bound. Note that the Sobolev space $\bar
{\mathbb H}(\alpha)=H^{\alpha+1/2}$ is contained in the Besov space
$B^{\alpha+1/2}_{22}([0,1])$, which itself is continuously imbedded
into the Besov space $B^{\alpha+1/2-1/2 +1/\bar r}_{\bar r 2}([0,1])=
B^{\alpha+1/\bar r}_{\bar r2}([0,1])$; cf. Remark \ref{besovmania}. We
conclude, for some constant $C'$, that $\|K_j(e^{f_1}) - e^{f_1}\|_r
\le C'\|f_1\|_{\bar{\mathbb H}(\alpha)} 2^{-j(\alpha+ 1/\bar r)}$ from
the approximation properties of wavelet projections on Besov spaces
(Definition \ref{besov}). This establishes the bound in the first part
of the lemma.

For the case $r=\infty$, note that, $f$ being bounded by $c$, the chain
rule gives that there exists $C(c,\alpha)$ such that
%
\begin{equation}\label{Tw}
\|e^{f}\|_{\alpha, \infty,1/2}\le C(c,\alpha)(\|f\|_{\alpha, \infty,1/2}+1).
\end{equation}
We conclude from a standard bias bound for wavelet projections that\break $\|
K_j(e^{f}) - e^{f}\|_\infty\le c(\|f\|_{\alpha, \infty, 1/2 }+1)
2^{-j\alpha} \sqrt j $ which, in view of $e^{-c}\le\int_0^1e^{f}\le
e^c$ gives the overall inequality.
\end{pf}

The choice $j=J_n$ with $2^{J_n} \sim n \varepsilon_n^2$, relevant in
Theorems \ref{general} and \ref{generalbd}, gives, for $p$ satisfying
the hypotheses of the previous proposition, the bounds
%
\begin{equation} \label{intbias0}
\|K_{J_n}(p)-p\|_r \le C \bigl((n \varepsilon_n^2)^{-\alpha}+\varepsilon
_n\bigr) \qquad\mbox{for }1 \le r \le2
\end{equation}
and
%
\begin{equation} \label{intbias}
\|K_{J_n}(p)-p\|_r \le C \bigl(\sqrt n \varepsilon_n (n \varepsilon
_n^2)^{-(\alpha+1/r)} + \varepsilon_n \bigr)\qquad\mbox{for }2<r<\infty
\end{equation}
as well as
%
\begin{equation} \label{intb}
\|K_{J_n}(p)-p\|_\infty\le C \sqrt n \varepsilon_n (n \varepsilon
_n^2)^{-\alpha} \sqrt{\log(n \varepsilon_n^2)}.
\end{equation}

The last auxiliary fact that we will require about $B_\alpha$ is a
small ball probability estimate, concretely an upper bound for the
concentration function~$\phi_w^{\bar B_\alpha}(\varepsilon)$ as
$\varepsilon$ approaches zero.
\begin{proposition}\label{smb}
Let $B_\alpha, \alpha\in\{n-1/2\dvtx n \in\mathbb N\}$ be integrated
Brownian motion, considered as a Gaussian vector taking values in the
Banach spa\-ce~$C([0,1])$, and let $w\in C^\alpha([0,1])$. Then, $\phi
_w^{\bar B_\alpha}(\varepsilon)=O(\varepsilon^{-1/\alpha})$, and the
same is true for~$\phi_w^{B_\alpha}$ if we further assume
$w^{(k)}(0)=0$, $k\le[\alpha]$.
\end{proposition}
\begin{pf}
Since $B_\alpha=W_{2\alpha}$ in \cite{LL98} and it also equals a
constant times $R_\alpha$ in~\cite{VV08}, this proposition simply
combines Theorem 2.1 in \cite{LL98} and Theorem~4.3 in \cite{VV08}.
\end{pf}

This result applies to the ``conditional'' concentration function: if $\|
w_0\|_\infty\le c/2$ and $\varepsilon\le c/2$, then
%
\begin{eqnarray}\label{cond}
&&
\Pr\{\|\bar B_\alpha-w_0\|_\infty<\varepsilon|\|\bar B_\alpha\|
_\infty\le c\}\nonumber\\
&&\qquad=
\frac{\Pr\{\|\bar B_\alpha-w_0\|_\infty<\varepsilon, \|\bar
B_\alpha\|_\infty\le c\}}{\Pr\{\|\bar B_\alpha\|_\infty\le c\}
}\\
&&\qquad=\frac{e^{-\phi_{w_0}^{\bar B_\alpha}(\varepsilon)}}{\Pr\{\|\bar
B_\alpha\|_\infty\le c\}}.\nonumber
\end{eqnarray}

We are now in a position to apply Theorems \ref{general} and \ref
{generalbd} to prove Proposition~\ref{fbmp}. To ease\vadjust{\goodbreak} notation define
$I(w)=e^{w}/\int_0^1e^{w(t)}\,dt, w\in C([0,1])$,
and record that, for $\|w\|_\infty\le c$,
%
\begin{equation} \label{ILips}
|I(w)| \le L (|w|+1),
\end{equation}
where $L$ depends only on $c$.

Set $w_0=\log p_0$, so that, since $\|w_0\|_\infty\le c/2$ and $p_0$ is
a density, hence $p_0=I(w_0)$,
Lemma 3.1 in \cite{VV08} gives that if $p=I(w)$ for $w=\bar B_\alpha
(\omega)$ for some $\omega\in\Omega$, and $\|w\|_\infty\le c$, then
$-P_0\log\frac{p}{p_0}\le R\|w-w_0\|_\infty^2$ and $P_0(\log\frac
{p}{p_0})^2\le R\|w-w_0\|_\infty^2$ for some $R<\infty$ (that depends
on $c$). Hence, for any $\varepsilon>0$ such that $R^{-1/2}\varepsilon< c/2$,
%
\begin{eqnarray}\label{KL}
&&\Pi\biggl\{p\in\mathcal{P}\dvtx -P_0\log\frac{p}{p_0}\le\varepsilon^2,
P_0\biggl(\log\frac{p}{p_0}\biggr)^2\le\varepsilon^2\biggr\}\nonumber\\[-8pt]\\[-8pt]
&&\qquad\ge\Pr\{\|\bar B_\alpha-w_0\|_\infty\le R^{-1/2}\varepsilon
|\|\bar B_\alpha\|_\infty\le c\}.\nonumber
\end{eqnarray}
Since $w_0$ is in $ C^{\alpha}([0,1])$, it follows from Proposition \ref
{smb} that
$\phi_{w_0}^{\bar B_\alpha}(\varepsilon)=O(\varepsilon^{-1/\alpha})$ as
$\varepsilon\to0$,
say, there exist $c_1$ large enough and $\varepsilon_1>0$ such that
\[
\phi_{w_0}^{\bar B_\alpha}(\varepsilon)\le c_1\varepsilon^{-1/\alpha}
\qquad\mbox{for all } \varepsilon\le\varepsilon_1.
\]
Then we have, for $\varepsilon_n = (c _1/n)^{\alpha/(2\alpha+1)}$, from
some $n$ on, both
\[
\phi_{w_0}^{\bar B_\alpha}(R^{-1/2}\varepsilon_n)\le c_1R^{1/(2\alpha
)}\varepsilon_n^{-1/\alpha} \quad\mbox{and}\quad  \phi_{w_0}^{\bar B_\alpha
}(\varepsilon_n)\le n\varepsilon_n^2.
\]
Hence, for these $n$, by (\ref{cond}),
%
\begin{equation}\label{TSMB}
\Pr\{\|\bar B_\alpha-w_0\|_\infty\le R^{-1/2}\varepsilon_n|
\|\bar B_\alpha\|_\infty\le c\}\ge e^{-Cn\varepsilon_n^2},
\end{equation}
where $C\,{=}\,c_1R^{1/(2\alpha)}$.
This proves condition (2) in Theorems \ref{general}, \ref{generalbd}
for these~$C, \varepsilon_n$.

To proceed with the verification of the conditions of Theorem \ref
{general}, take
$\mathcal P_n = \{I(w)\dvtx w \in\{M_n \bar{\mathbb H}_1(\alpha) +
\varepsilon_n B^{1}\}\}$ if $r<\infty$ and
$\mathcal P_n=\{I(w)\dvtx \|w\|_{\alpha,\infty,1/2} \le\sqrt{(C+4)/L_\alpha
} \sqrt n \varepsilon_n\}$ if $r=\infty$, and note that condition (1)
in Theorem \ref{general} is satisfied for these choices in view of
Propositions \ref{borelll0} and \ref{borell}; see (\ref{borell0}) and~(\ref{Borellll}). The bias condition is satisfied for the above choice
of $\varepsilon_n$, $\gamma_n=1$ if $r<\infty$ and $\gamma_n= \sqrt
{\log n}$ if $r=\infty$, in view of Lemma \ref{BIAS}; cf. also (\ref
{intbias0}), (\ref{intbias}), (\ref{intb}). Finally\vadjust{\goodbreak} the additional
restrictions on $\varepsilon_n$ in Theorems \ref{general} and \ref
{generalbd} are also satisfied, unless $\alpha=1/2, r=\infty$. In this
case the rate of contraction $\delta_n$ exceeds a~constant multiple
times $\sqrt{\log n}$, so that the result follows trivially from the
fact that the prior is supported in a sup-norm bounded set.

\subsection{\texorpdfstring{Proof of Theorem \protect\ref{whitenoise1}}{Proof of Theorem 1}}

Observing $Y^{(n)}$ is equivalent to observing its action, on the basis,
%
\begin{eqnarray}\label{y1}
y_k&=&\int_0^1\phi_k(t)\,dY^{(n)}(t)=\langle f,\phi_k\rangle+\frac
{1}{\sqrt{n}}\int_0^1\phi_k(t)\,dB(t)\nonumber\\[-9pt]\\[-9pt]
:\!&=&\theta_k+\frac{1}{\sqrt{n}}g_k,\qquad k=0,\ldots, N-1,\nonumber\\[-2pt]
%
\label{y2}
y_{\ell k}&=&\int_0^1\psi_{\ell k}(t)\,dY^{(n)}(t)\nonumber\\[-2pt]
&=&\langle f,\psi_{\ell
k}\rangle+\frac{1}{\sqrt{n}}\int_0^1\psi_{\ell
k}(t)\,dB(t)\\[-2pt]
:\!&=&\theta_{\ell k}+\frac{1}{\sqrt{n}}g_{\ell k}, \qquad k=0,\ldots,2^{\ell
}-1,\ell\ge J_0,\nonumber
\end{eqnarray}
with the variables $g_k$, $g_{\ell k}$ all i.i.d. $N(0,1)$. The
observed process, still denoted by $Y^{(n)}$, can thus be viewed as a
random element $Y^{(n)}=(y_k,y_{\ell k})^t$ of~$\ell_2$, where $y_k$ is
$N(\theta_k,1/n)$, and $y_{\ell k}$ is $N(\theta_{\ell k},1/n)$, all
independent. Likewise\vspace*{1pt} the function $f_0$ to be estimated becomes the
vector $\theta_0=(\theta_{k}^0,\theta_{\ell k}^0)^t$ of the
coefficients of its wavelet expansion, that is, $\theta_k^0=\langle
f_0,\phi_k\rangle$ and $\theta_{\ell k}^0=\langle f_0,\phi_{\ell
k}\rangle$, and any prior $\Pi$ on $L_2$ maps onto a prior, still
denoted by $\Pi$, on the parameter space $\theta=(\theta_k,\theta_{\ell
k})^t\in\ell_2$.

The posterior $\Pi(\cdot|Y^{(n)})$ is then the law of $\theta$ given
the observed process~$Y^{(n)}$. Standard results on Gaussian measures
on $\ell^2$ imply that if the prior $\Pi$ on $\ell_2$ is a centered
Gaussian vector of trace class covariance $\Sigma$, then the posterior
probability law given $Y^{(n)}$, $\hat\Pi_n^Y=\hat\Pi^{Y^{(n)}}$, is
also Gaussian, with mean $\hat\theta(Y)=E_{\Pi}(\theta|Y^{(n)})=\Sigma
(\Sigma+I/n)^{-1}Y^{(n)}=\Sigma(\Sigma+I/n)^{-1}(y_k;y_{\ell k})^t$ and
with\vspace*{1pt} covariance $\Sigma|Y^{(n)}=\Sigma(n\Sigma+I)^{-1};$ see, for
example, Theorem 3.2 in~\cite{Z00}. We will drop the superindex $(n)$
from the processes $Y^{(n)}$ and $Y_0^{(n)}$ from now on to expedite
notation.\vspace*{-0.5pt}

The posterior $\hat\Pi_n^Y$ gives rise to a Gaussian measure on
$L_2([0,1])$ by simply ``undoing'' the isometry, that is, by taking the
law\vspace*{1pt} of the random wavelet series in $L^2([0,1])$ with coefficients
drawn from $\hat\Pi_n^Y$ equal to
\begin{eqnarray*}
X&=&\sum_{k=0}^{N-1}\biggl[\frac{1}{1+1/n}y_k+\biggl(\frac{1}{n+1}\biggr)^{1/2}
\bar g_k\biggr]\phi_k\\[-2pt]
&&{}+\sum_{\ell=J_0}^{\infty}\sum_{k=0}^{2^{\ell}-1}\biggl[\frac{\mu_{\ell
}}{\mu_{\ell}+1/n}y_{\ell k}+ \biggl(\frac{\mu_\ell}{n\mu_\ell+1}
\biggr)^{1/2}\bar g_{\ell k}\biggr]\psi_{\ell k}\\[-2pt]
&=& E_{\Pi_n}(f|Y)+\sum_{k=0}^N\biggl(\frac{1}{n+1}\biggr)^{1/2}
\phi_k\bar g_k\\[-3pt]
&&{}+\sum_{\ell=J_0}^\infty\sum_{k=0}^{2^{\ell}-1}\biggl(\frac
{\mu_\ell}{n\mu_\ell+1}\biggr)^{1/2}\psi_{\ell k}\bar g_{\ell
k},\vspace*{-2pt}
\end{eqnarray*}
where the $\bar g$ variables are i.i.d. $N(0,1)$, and $y_k$, $y_{\ell
k}$ are, as defined above, the integrals of the wavelet basis functions
with respect to $dY(t)$. Under $dY_0(t)=f_0(t)\,dt+dB(t)/\sqrt{n}$, we
have $y_k=\langle f_0,\phi_k\rangle+g_k/\sqrt{n}$, $y_{\ell k}=\langle
f_0,\phi_{\ell k}\rangle+g_{\ell k}/\sqrt{n}$, where the $g_k, g_{\ell
k}$ are again i.i.d. $N(0,1)$, independent of the variables $\bar g$.
So, the posterior given $Y_0$ integrates the $\bar g$ variables, and
$E_{Y_0}$ integrates the $g$ variables, and we have
%
\begin{eqnarray}\label{exp}
&&E_{Y_0}\hat\Pi_n^{Y_0}\{\|f-f_0\|_\infty>M\varepsilon_n\}
\nonumber\\[-3pt]
&&\qquad=\Pr\Biggl\{\Biggl\|\sum_{k=0}^{N-1}\biggl[\frac{-1/n}{1+1/n}\langle
f_0,\phi_k\rangle\nonumber\\[-5pt]
&&\qquad\quad\hspace*{44.3pt}{}+\frac{1}{\sqrt{n}(1+1/n)}g_k+\biggl(\frac{1}{n+1}
\biggr)^{1/2}\bar g_k\biggr]\phi_k\nonumber\\[-5pt]
&&\qquad\quad\hphantom{\Pr\Biggl\{\Biggl\|}
{}+\sum_{\ell=J_0}^{\infty}\sum_{k=0}^{2^\ell-1}\biggl[\frac{-1/n}{\mu
_{\ell}+1/n}\langle f_0,\psi_{\ell k}\rangle\\[-5pt]
&&\qquad\quad\hspace*{77pt}{}+\frac{\mu_{\ell}}{\sqrt
{n}(\mu_{\ell}+1/n)}g_{\ell k}\nonumber\\[-5pt]
&&\hspace*{106pt}\qquad\quad{}+
\biggl(\frac{\mu_\ell}{n\mu_\ell+1}\biggr)^{1/2}\bar g_{\ell k}\biggr]\psi
_{\ell k} \Biggr\|_\infty>M\varepsilon_n\Biggr\}\nonumber\\[-5pt]
&&\qquad= \Pr\bigl\{\bigl\|E_{Y_0^{(n)}}\bigl(E_{\Pi_n}
(f|Y_0)-f_0\bigr)+G\bigr\|
_\infty>M\varepsilon_n\bigr\},\nonumber\vspace*{-2pt}
\end{eqnarray}
where $G$ is the centered Gaussian process
\begin{eqnarray*}
G(t)&=&\sum_{k=0}^{N-1}\biggl[\frac{1}{\sqrt{n}(1+1/n)}g_k+\biggl(\frac
{1}{n+1}\biggr)^{1/2}\bar g_k\biggr]\phi_k(t)\\[-3pt]
&&{}+\sum_{\ell=J_0}^{\infty}\sum_{k=0}^{2^\ell-1}\biggl[\frac{\mu
_{\ell}}{\sqrt{n}(\mu_{\ell}+1/n)}g_{\ell k}+
\biggl(\frac{\mu_\ell}{n\mu_\ell+1}\biggr)^{1/2}\bar g_{\ell k}
\biggr]\psi_{\ell k}(t)\vspace*{-2pt}
\end{eqnarray*}
and
\begin{eqnarray*}
E_{Y_0}\bigl(E_{\Pi_n}(f|Y_0)-f_0\bigr)&=&\sum_{k=0}^{N-1}\frac
{-1/n}{1+1/n}\langle f_0,\phi_k\rangle\phi_k\\[-3pt]
&&{} +\sum_{\ell=J_0}^\infty\sum_{k=0}^{2^\ell-1}\frac{-1/n}{\mu_{\ell
}+1/n}\langle f_0,\psi_{\ell k}\rangle\psi_{\ell k}.\vspace*{-2pt}
\end{eqnarray*}
It suffices to prove the theorem for $r =\infty$. We will apply
Borell's \cite{B75} inequality (a\vadjust{\goodbreak} consequence thereof, in fact,
equation (3.2) in \cite{LT91}, page 57) to the probability in (\ref
{exp}), and for this we need to estimate $\|E(E_{\Pi_n}(f|Y_0)-f_0)\|
_\infty$, $E\|G\|_\infty$ and $\|E(G^2(\cdot))\|_\infty$.

Choose $J_n\ge J_0$ such that $2^{J_n} \simeq(n/\log n)^{1/(2\alpha
+1)}$. Since $f_0 \in\mathcal C^{\alpha}([0,1])$ and $\|{\sum
_k}|\psi_{\ell k}|\|_\infty\le C2^{\ell/2}$, we obtain
\[
\Biggl\|\sum_{k=0}^{N-1}\frac{-1/n}{1+1/n}\langle f_0,\phi_k\rangle\phi
_k\Biggr\|_\infty\le \Biggl\|\sum_{k=0}^{N-1}|\phi_k|\Biggr\|_\infty
\frac{C}{n+1}\le\frac{ C_1}{n}\vspace*{-2pt}
\]
and
\begin{eqnarray*}
\Biggl\|\sum_{\ell=J_0}^\infty\sum_{k=0}^{2^\ell-1}\frac{-1/n}{\mu_{\ell
}+1/n}\langle f_0,\psi_{\ell,k}\rangle\psi_{\ell,k}\Biggr\|_\infty
&\le&
\sum_{\ell=J_0}^\infty\Biggl\|\sum_{k=0}^{2^\ell-1}|\psi_{\ell
k}|\Biggr\|_\infty\frac{C2^{-\ell(\alpha+1/2)}}{n\mu_\ell+1} \\[-2pt]
&\le& C'\Biggl(\sum_{\ell=J_0}^{J_n}\frac{2^{-\ell\alpha}}{n \mu_\ell}
+\sum_{\ell=J_n+1}^\infty2^{-\ell\alpha}\Biggr)\\[-2pt]
& \le& C_2\biggl(\frac{\log n}{n}\biggr)^{\alpha/(2\alpha+1)},\vspace*{-2pt}
\end{eqnarray*}
where $C_1$ and $C_2$ depend only on the wavelet basis, $\alpha$ and $\|
f_0\|_{\alpha, \infty}$.
Collecting the last two sets of inequalities yields the bound
%
\begin{equation}\label{E1}
\bigl\|E_{Y_0}\bigl(E_{\Pi_n}(f|Y_0)-f_0\bigr)\bigr\|_\infty\le\bar C_1
\biggl(\frac{\log n}{n}\biggr)^{\alpha/(2\alpha+1)}\vspace*{-2pt}
\end{equation}
for some $\bar C_1<\infty$. To bound $E\|G\|_\infty$, recall that for
any sequence of centered normal random variables $Z_j$,
%
\begin{equation}\label{bore}
E\max_{1\le j\le N}|Z_j|\le C\sqrt{\log N}\max_{j\le N}(EZ_j^2)^{1/2},\vspace*{-2pt}
\end{equation}
where $C$ is a universal constant. Therefore, from the definitions of
$J_n, \mu_\ell$,
\begin{eqnarray*}
&& E\biggl\|\sum_k\biggl[\frac{1}{\sqrt{n}(1+1/n)}g_k+\biggl(\frac
{1}{n+1}\biggr)^{1/2}\bar g_k\biggr]\phi_k\biggr\|_\infty\\[-2pt]
&&\qquad \le\biggl\|\sum_k|\phi_k|\biggr\|_\infty\biggl(\frac
{1}{n(1+1/n)^2}+\frac{1}{n+1}\biggr)^{1/2}E\max_k|g_k| \\[-2pt]
&&\qquad=O\biggl(\frac
{1}{\sqrt n}\biggr)\vspace*{-2pt}
\end{eqnarray*}
and, using $\mu_\ell\lesssim n^{-1}$ for $\ell\ge J_n$,
\begin{eqnarray*}
&&E\Biggl\|\sum_{\ell=J_0}^\infty\sum_{k=0}^{2^\ell-1}\biggl[\frac{\mu
_{\ell}}{\sqrt{n}(\mu_{\ell}+1/n)}g_{\ell k}+
\biggl(\frac{\mu_\ell}{n\mu_\ell+1}\biggr)^{1/2}\bar g_{\ell k}
\biggr]\psi_{\ell k}\Biggr\|_\infty\\[-2pt]
&&\qquad \le C'\sum_{\ell=J_0}^\infty2^{\ell/2}E \max_{k \le2^\ell}|g_{\ell
k}|\biggl(\frac{\mu_\ell^2}{n(\mu_\ell+1/n)^2}+\frac{\mu_\ell}{n\mu_\ell
+1}\biggr)^{1/2} \\[-2pt]
&&\qquad\le C''\sum_{\ell=J_0}^\infty(\ell2^\ell)^{1/2}\biggl(\frac{\mu_\ell
^2}{n(\mu_\ell+1/n)^2}+\frac{\mu_\ell}{n\mu_\ell+1}\biggr)^{1/2} \\[-2pt]
&&\qquad\le C''\Biggl(2\sum_{\ell=J_0}^{J_n}\sqrt{\frac{2^\ell\ell}{n}}+\sum
_{\ell> J_n}\sqrt{2^\ell\ell n}\mu_\ell+\sum_{\ell> J_n}\sqrt{2^\ell
\ell\mu_\ell}\Biggr)\\[-2pt]
&&\qquad\le C'''\Biggl(\sqrt{\frac{2^{J_n}J_n}{n}}+2^{-J_n \alpha}\Biggr) \le
D\biggl(\frac{\log n}{n}\biggr)^{\alpha/(2\alpha+1)}.
\end{eqnarray*}
Conclude
%
\begin{equation}\label{gnorm}
E\|G\|_\infty\le\bar C_2\biggl(\frac{\log n}{n}\biggr)^{\alpha/(2\alpha+1)}
\end{equation}
for some $\bar C_2<\infty$.
Finally,
%
\begin{eqnarray}\label{eg2}
EG^2(t)&=&\sum_{k=0}^{N-1}\biggl(\frac{1}{n(1+1/n)^2}+\frac{1}{n+1}
\biggr)\phi_k^2(t) \nonumber\\
&&{} +\sum_{\ell=J_0}^\infty\sum_{k=0}^{2^\ell-1}\biggl(\frac{\mu_{\ell
}^2}{n(\mu_{\ell}+1/n)^2}+\frac{\mu_\ell}{n\mu_\ell+1}\biggr)\psi_{\ell
k}^2(t)\\
&\le& C\biggl(\frac{1}{n}+ \frac{2^{J_n}}{n}+2^{-J_n(2\alpha+1)}\biggr)
\le C_3\frac{2^{J_n}}{n}.\nonumber
\end{eqnarray}
So, setting $\varepsilon_n=(n/\log n)^{-\alpha/(2\alpha
+1)}$, the estimates (\ref{E1}), (\ref{gnorm}) and (\ref{eg2}) together
with inequality (3.2) on page 57 of \cite{LT91}, give
%
\begin{eqnarray}\label{exp2}\quad
&&\Pr\bigl\{\|E_{Y_0}\bigl(E_{\Pi_n}(f|Y_0)-f_0\bigr)+G\|_\infty
>M\varepsilon_n\bigr\}\nonumber\\
&&\qquad\le\Pr\{\|G\|_\infty-E\|G\|_\infty>M\varepsilon_n-\bigl\|E
\bigl(E_{\Pi_n}(f|Y_0)-f_0\bigr)\bigr\|_\infty-E\|G\|_\infty\}
\nonumber\\[-8pt]\\[-8pt]
&&\qquad\le\Pr\{\|G\|_\infty-E\|G\|_\infty>(M-\bar C_1-\bar
C_2)\varepsilon_n\}\nonumber\\
&&\qquad\le\exp\biggl(-\frac{(M-\bar C_1-\bar C_2)^2\varepsilon
_n^2}{C_3^22^{J_n}/n}\biggr).\nonumber
\end{eqnarray}
Collecting (\ref{exp}) and (\ref{exp2}) and taking into account that
$\varepsilon_n^2\simeq2^{J_n}J_n/n$ completes the proof.

\section*{Acknowledgments}

E. Gin\'{e}'s work was carried out during a sabbatical leave at the
M.I.T. Mathematics Department and on a visit at the University of Cambridge
Statistical Laboratory, and he is grateful for the hospitality of these
institutions.
R. Nickl would like
to thank the Caf\'{e}s Br\"{a}unerhof and Florianihof
in Vienna for their continued hospitality.

We are further grateful to Ismael Castillo, Vladimir Koltchinskii, Natesh
Pillai, Catia Scricciolo and Aad van der Vaart for valuable
conversations about the subject of this article. We also thank the
Associate Editor and two referees for substantial and influential
reports on preliminary versions of this manuscript.


\begin{supplement} 
\stitle{Supplement to ``Rates of contraction for posterior
distributions in $\bolds{L^r}$-metrics, $\bolds{1 \le r \le \infty}$''}
\slink[doi]{10.1214/11-AOS924SUPP} 
\sdatatype{.pdf}
\sfilename{aos924\_supp.pdf}
\sdescription{This supplement contains a detailed proof
of Lemma \ref{momr} and an expanded proof of Proposition \ref{dir}
from the mentioned article.}
\end{supplement}


\printaddresses

\end{document}